\documentclass[reqno]{amsart}
\usepackage{a4wide}
\usepackage[colorlinks=true]{hyperref}

\theoremstyle{plain}
\newtheorem{thm}{Theorem}[section]
\newtheorem{proposition}[thm]{Proposition}

\newtheorem{corollary}[thm]{Corollary}

\theoremstyle{definition}
\newtheorem{definition}[thm]{Definition}
\newtheorem{example}[thm]{Example}

\theoremstyle{remark}
\newtheorem{remark}[thm]{Remark}

\newcommand\Mink{\ensuremath{\mathbb{R}^{(2, 1)}}}
\newcommand\nequiv{\mathrel/ \hspace{-6pt}\equiv}

\begin{document}
\title{Coisotropic variational problems}
\author{Emilio Musso}
\author{James D.E.\ Grant}
\email{\href{mailto:musso@univaq.it}{musso@univaq.it}; 
\href{http://www.emiliomusso.com}{http://www.emiliomusso.com}}
\email{\href{mailto:grant@dm.univaq.it}{grant@dm.univaq.it}} 
\address{Dipartimento di Matematica Pura ed Applicata\\
Universit{\`a} di L'Aquila\\
Via Vetoio\\
67100 L'Aquila\\
Italy.}
\begin{abstract}
In this article we study constrained variational problems in
one independent variable defined on the space of integral curves of a
Frenet system in a homogeneous space $G/H$.
We prove that if the Lagrangian is $G$-invariant and coisotropic then the
extremal curves can be found by quadratures.
Our proof is constructive and relies on the
reduction theory for coisotropic optimal control problems.
This gives a unified explanation of the integrability
of several classical variational problems such as the
total squared curvature functional, the projective,
conformal and pseudo-conformal arc-length functionals,
the Delaunay and the Poincar{\'e} variational problems.
\end{abstract}
\subjclass[2000]{Primary $58$A$30$, $53$D$20$; Secondary $58$A$10$, $37$K$10$.}
\thanks{Partially supported by MURST project
\lq\lq Propriet{\`a} Geometriche delle Variet{\`a} Reali e Complesse\rq\rq, by
GNSAGA and by European Contract Human Potential Programme,
Research Training Network HPRN-CT-$2000$-$00101$ (EDGE). 
Part of this work was done while the first author was
visiting the Laboratoire de Math{\'e}matiques J.A.Dieudonn\'{e},
Universit{\'e} de Nice Sophia-Antipolis.
He would like to thank this institution for their support and hospitality.}

\date{15 July, 2003}
\maketitle
\thispagestyle{empty}

\section{Introduction}
The present paper is an outcome of our attempt to understand the
general mechanisms underlying the integrability of constrained
variational problems for curves of constant type in homogeneous
spaces \cite{Green,Griffiths2,Jensen1,Jensen2}. The Pfaffian
differential systems arising from curves of constant type lead to
the notion of \emph{generalized Frenet system\/} for curves of a
homogeneous space $G/H$. Roughly speaking, a generalized Frenet
system of order $k$ on $G/H$ is a $G$-invariant submanifold
$S \subset J^k(\mathbb{R}, G/H)$ of the jet space $J^k(\mathbb{R}, G/H)$,
which may be linearized by a left-invariant affine sub-bundle of $T(G)$. From
the geometrical viewpoint the integral curves of such systems are
canonical lifts of curves of constant type on $G/H$. The most
elementary example is the classical Frenet-Serret differential
system for generic curves in Euclidean space. We then
consider a $G$-invariant Lagrangian and we investigate the
corresponding Euler-Lagrange system. The general construction of
the momentum space and of the Euler-Lagrange system of a
constrained variational problem in one independent variable is due
to Griffiths (we refer to \cite{BryantGriffiths,Griffiths1,Hsu}
as the standard references on the subject and to \cite{Cartan1} as
the original source of inspiration of the approach developed by
Griffiths). We adhere to the terminology introduced in
\cite{Griffiths1,BryantGriffiths} and say that a Lagrangian $L$
is \emph{non-degenerate\/} if the \emph{momentum space\/} $Y$ is
odd-dimensional and if the canonical $2$-form on $Y$ has maximal
rank. We prove that if the Lagrangian $L$ is $G$-invariant and
\emph{coisotropic\/} (see Definition~\ref{def:coisotropic})
then the extremal curves of the variational problem can be found
by quadratures. The proof relies on the reduction theory of
Hamiltonian systems with symmetries (see
\cite{ArnoldGivental,FomenkoTrofimov,GuilleminSternberg1,GuilleminSternberg2}
for the standard theory in the symplectic category and
\cite{Albert,BryantGriffiths,OrtegaRatiu,Willett} for
generalizations to contact geometry, time-dependent Hamiltonian
systems and Poisson manifolds). One of the ingredients of the
proof is a concrete geometric description of the Marsden-Weinstein
reduced spaces in terms of the \emph{phase portraits\/} of the
system. This procedure is constructive and applies to several
concrete examples (see refs.
\cite{BryantGriffiths,Cartan2,Griffiths1,Hsu,Musso1,Musso2,Musso3,MussoNicolodi}).

The paper is organized as follows. In the
next section we recall the basic definitions and properties of
linear control systems on Lie groups and Frenet systems
of curves in homogeneous spaces. In Section~\ref{sec:variational}, 
we examine variational problems
defined by invariant Lagrangians for linear control systems on Lie
groups. From a geometrical viewpoint we deal with $k$-th order
variational problems for curves of constant type in a homogeneous
space that depend on the generalized curvatures. Since all the
derived systems have constant rank, the extremal curves of
the variational problem are the projections of the integral curves
of the Euler-Lagrange system. Therefore, we focus our attention on the
momentum space and on the Euler-Lagrange system. First we
investigate the geometry of the momentum space $Y$ of a 
\emph{regular invariant Lagrangian\/} 
of a linear control system of a Lie group
$G$. We show that $Y$ is of the form $G \times \mathcal{F} $, 
where $\mathcal{F}$ is an immersed submanifold of $\mathfrak{g} \times \mathfrak{g}^*$ 
(we call $\mathcal{F}$ the \emph{phase space\/} of the variational problem). 
Next we study non-degenerate Lagrangians. 
We prove that if $L$ is non-degenerate, 
then the phase space $\mathcal{F}$ can be realized as a submanifold of $\mathfrak{g}^*$. 
We define the linearized phase portraits and the Legendre transform and 
analyze the structure of the characteristic vector field of a non-degenerate Lagrangian.
In Section~\ref{sec:coisotropic} we study \emph{coisotropic Lagrangians}. 
We prove that the integral curves of the characteristic vector field 
passing through a point of the bifurcation set 
are orbits of one-parameter subgroups of the symmetry group $G$. 
Therefore, from this point on, we focus our attention on the regular part $Y_r$ of the momentum space. 
We show that $Y_r$ is of the form $G \times \mathcal{F}_r$, 
where $\mathcal{F}_r$ is an open subset of the phase space. 
We prove that $\mathcal{F}_r$ intersects the coadjoint orbits, $\mathcal{O}(\mu)$, of $G$ transversally 
and that $\mathcal{P}_r(\mu) = \mathcal{F}_r \cap \mathcal{O}(\mu)$ are smooth curves
(referred to as the \emph{phase portraits\/}). 
Subsequently we introduce the moment map $J: Y_r \to \mathfrak{g}^*$ 
and prove that the Marsden-Weinstein reduction $J^{-1}(\mu)/G_{\mu}$ can be naturally
identified with the phase portrait $\mathcal{P}_r (\mu)$. 
We also show that every $\mu \in J(Y_r)$ is a regular element of $\mathfrak{g}^*$ 
which implies that the isotropy subgroups $G_{\mu}$ are Abelian, 
for every $\mu \in J(Y_r)$. 
We then examine more closely the phase flow $\phi$ and the characteristic vector field $\xi$. 
We prove that if the Lie algebra $\mathfrak{g}$ possesses a non-degenerate 
$Ad$-invariant inner product then the differential equation fulfilled by the phase flow
can be written in Lax form. From the Noether conservation theorem
we know that the characteristic vector field $\xi$ is tangent to
the fibers $J^{-1}(\mu)$ of the moment map. We define a canonical
connection form $\theta^{\mu}$ on the Marsden-Weinstein fibrations
$J_{r}^{-1}(\mu) \to \mathcal{P}(\mu)$ whose horizontal curves are the
integral curves of the characteristic vector field. Since the base
is $1$-dimensional and the structure group $G_{\mu}$ is Abelian, 
the horizontal curves can be found by a single quadrature.
This shows that the extremal curves of an invariant coisotropic
Lagrangian are integrable by quadratures. 
As a byproduct, we prove that if the canonical connection $\theta_{\mu}$ is complete, 
which is generically the case when $\mu$ is a regular value of the moment map, 
then the connected components of $J^{-1}_r(\mu)$ are Euclidean cylinders 
and the characteristic vector field $\xi$ can be linearized on $J^{-1}_r(\mu)$. 
We would like to stress that the connection form $\theta^{\mu}$ 
can be constructed explicitly from the data of the problem, 
so that the integration process can be performed in a completely explicit way. 

Finally, in two appendices, we summarise the background material that 
we use from the theory of Pfaffian differential systems and constrained 
variational problems in one independent variable. 

Throughout the paper, we demonstrate how our general results apply to the specific example 
of isotropic curves in $\Mink$. We show how to derive the Frenet system for such curves, 
and show that the variational problem is coisotropic if we take the Lagrangian to be a linear function 
of the curvature. We prove that the phase portraits may be parametrised in terms of elliptic 
functions, and construct the sections of the Marsden-Weinstein fibration required to reduce 
the integration to quadratures. Other concrete geometrical examples where the general scheme 
described in this paper are implemented may be found in 
\cite{BryantGriffiths,Cartan1,Musso1,Musso2,Musso3,MussoNicolodi}. 
In all of these cases the generic phase portrait is an elliptic curve, 
so that the extremal curves can be integrated in terms of elliptic functions and elliptic integrals.

\section{Linear control systems on Lie groups and Frenet systems in homogeneous spaces}
\subsection{Linear control systems on Lie groups}
Let $G$ be a Lie group with Lie algebra $\mathfrak{g}$.
The natural pairing between $\mathfrak{g}$ and $\mathfrak{g}^*$ will be denoted by
$(\eta, V) \in \mathfrak{g}^* \times \mathfrak{g} \rightarrow \langle \eta; V \rangle \in \mathbb{R}$.
We let $\Theta \in \mathfrak{g}^* \otimes \mathfrak{g}$ be the Maurer-Cartan form of $G$.
If we fix a basis $(e_0, \dots, e_n)$ of $\mathfrak{g}$, then $\Theta = \theta^J \otimes e_J$,
where $(\theta^0, \dots, \theta^n)$ is the basis of $\mathfrak{g}^*$ dual to $(e_0, \dots, e_n)$.

\begin{definition}
Let $\mathbb{A} \subset \mathfrak{g}$ denote an affine subspace of $\mathfrak{g}$ of the form 
$P + \mathfrak{a} = \left\{ P + A: A \in \mathfrak{a} \right\}$, 
where $P \in \mathfrak{g}$ and $\mathfrak{a} \in Gr_h(\mathfrak{g})$ with $P \notin \mathfrak{a}$. 
The set of such affine subspaces of $\mathfrak{g}$ will be denoted $P^h(\mathfrak{g})$. 
We call $M:= G \times \mathbb{A}$ the \emph{configuration space\/} of the affine subspace $\mathbb{A}$, 
and denote by $\pi_G: M \to G$ and $\pi_{\mathbb{A}}: M \to \mathbb{A}$ 
the natural projections onto the two factors. 
\end{definition}

We now fix a left-invariant form $\omega \in \mathfrak{g}^*$ such that
$\langle \omega; P \rangle = 1$ and $\omega \in \mathfrak{a}^{\perp}$ 
(i.e. $\langle \omega; A \rangle = 0$, for all $A \in \mathfrak{a}$). 
We may fix a basis $(e_0, \dots, e_n)$ of $\mathfrak{g}$ such that
\[
P = e_0, \quad \mathfrak{a} = \mathrm{span}(e_1, \dots, e_h)
\]
and we let $\theta^0, \dots, \theta^n$ be the components of the Maurer-Cartan
form with respect to $(e_0, \dots, e_n)$. Such a basis may be chosen
so that $\theta^0 = \omega$. Using the projection $\pi_G$, 
we may pull-back the differential $1$-forms $\omega, \theta^1, \dots, \theta^n$ to $M$, 
to define a set of $1$-forms on $M$ which, by the standard abuse of notation, we 
again denote by $\omega, \theta^1, \dots, \theta^n \in \Omega^1(M)$
\footnote{We will generally follow the usual practice 
in the method of moving frames and omit the pull-back signs 
to simplify notation. This should cause no confusion as 
we will clearly specify the manifolds that we are working on.}. 
Let $k^1, \dots, k^h$ denote the affine coordinates on $\mathbb{A}$ defined by the affine frame $(P, e_1, \dots, e_h)$. 
We then define the $1$-forms
\[
\eta^j := 
\begin{cases} 
\theta^j - k^j \omega &j = 1, \dots, h,\\
\theta^j &j = h+1, \dots, n.
\end{cases}
\]
We then define the Pfaffian differential system $(\mathcal{A}, \omega)$ on $M$ 
to be the Pfaffian differential ideal generated by the $1$-forms $\{ \eta^j: j = 1, \dots n \}$ 
with the independence condition given by $\omega$
\footnote{More invariantly, if we fix $\omega \in \mathfrak{g}^*$ with 
$\omega \in \mathfrak{a}^{\perp}$ and $\langle \omega; P \rangle = 1$, 
then we define the $\mathfrak{g}$-valued $1$-form 
$\widehat{\Theta} \in \Omega^1 \left( M, \mathfrak{g} \right)$ by the formula 
$\widehat{\Theta}|_{(g, Q)} := \pi_G^*(\Theta - Q \omega)|_{(g, Q)}$. 
$\mathcal{A}$ is then the differential ideal generated by
$\{ \langle \mu; \widehat{\Theta} \rangle: \mu \in \mathfrak{g}^* \}$.}.

\begin{definition}
$(\mathcal{A}, \omega)$ is the \emph{linear control system associated to $\mathbb{A} \in P^h(\mathfrak{g})$}.
\end{definition}

Note that the ideal $\mathcal{A}$ has constant rank, being generated by a rank $n$ sub-bundle, $Z \subset T^*(M)$.
The sub-bundle $Z$ is of the form $G \times \mathcal{Z}$, where
\[
\mathcal{Z} = \{(Q, \eta) \in \mathbb{A} \times \mathfrak{g}^*: \langle \eta; Q \rangle = 0 \} 
\subset \mathbb{A} \times \mathfrak{g}^*,
\]
and where the embedding of $Z$ as a sub-bundle of $T^*(M)$ is given by
\[
(g, Q, \eta) \in G \times \mathcal{Z} \to \left. \pi_G^*(\eta) \right|_{(g, Q)} \in T^*_{(g, Q)}M .
\]

We may use the left-invariant trivialization to identify $T(G)$ and $G \times \mathfrak{g}$.
The tangent space to $M$ at $(g, Q)$ is then identified with $\mathfrak{g} \oplus \mathfrak{a}$.
With this identification at hand, the integral elements of $(\mathcal{A}, \omega)$
at $(g, Q)$ are the $1$-dimensional subspaces of $\mathfrak{g} \oplus \mathfrak{a}$ of the form
$(Q, v)$, where $v \in \mathfrak{a}$.

A smooth curve $\gamma = (\alpha, \beta): (a, b) \to M$, 
where $(a, b) \subseteq \mathbb{R}$, 
is a parametrized integral curve of the control system $(\mathcal{A}, \omega)$ if and only if
$\alpha: (a, b) \to G$ is a solution of the linear system $\alpha(t)^{-1} \alpha'(t) = \beta(t)$.
Thus, as a control system, the points of the affine space $\mathbb{A}$ play the role of the inputs.
Note that if we assign a smooth map $\beta: (-\epsilon, \epsilon) \to \mathbb{A}$ and a point $g_0 \in G$, 
then there exists a unique integral curve of the control system, $\gamma = (\alpha, \beta)$, satisfying
the initial condition $\alpha(0) = g_0$.

Consider the linear subspaces $\mathfrak{a}_k \subset \mathfrak{g}$ defined recursively by 
\[
\mathfrak{a}_1 = \mathfrak{a} + \mathrm{span}(P), \quad
\mathfrak{a}_2 = \mathfrak{a}_1 + [\mathfrak{a}_1, \mathfrak{a}_1], \quad \cdots, \quad
\mathfrak{a}_k = \mathfrak{a}_{k-1} + [\mathfrak{a}_{k-1}, \mathfrak{a}_{k-1}].
\]
The smallest integer $N$ such that $\mathfrak{a}_N = \mathfrak{a}_{N+1}$ is called the \emph{derived length\/} of $\mathbb{A}$.
Note that $\mathcal{Z}_s:= \mathbb{A} \times \mathfrak{a}_s^{\perp}$ is contained in $\mathcal{Z}$, for $s = 1, \dots, N$.
We set $Z_s = M \times \mathcal{Z}_s$, $s = 1, \dots, N$ and we consider the sequence of sub-bundles
\[
Z_N \subset Z_{N-1} \subset \dots \subset Z_1 \subset Z.
\]
If we denote by $\mathcal{A}_s$ the Pfaffian differential ideal generated by $Z_s$,
then
\[
\mathcal{A}_N \subset \mathcal{A}_{N-1} \subset \dots \subset \mathcal{A}_1 \subset \mathcal{A}
\]
is the derived flag of the control system (see ref.~\cite{BCGG,Griffiths1} for more details about derived flags).
We have thus proved the following:

\begin{proposition}
All the derived systems of a linear control system on a Lie group $G$ have constant rank.
\end{proposition}

\subsection{Frenet systems in homogeneous spaces}
Let $H \subset G$ be a closed Lie subgroup and consider the
homogeneous space $G/H$. The left-action of $G$ on $G/H$
induces an action of $G$ on the jet space $J^k(\mathbb{R}, G/H)$, called
the \emph{$k$-th prolongation\/} of the action of G on $G/H$.

\begin{definition}
A \emph{differential relation (in one independent variable) of order $k$\/}
on $G/H$ is a submanifold $S$ of $J^k(\mathbb{R}, G/H)$ such that
$\left. dt \right|_S$ is nowhere vanishing. 
We define the Pfaffian differential system $(\mathcal{I}, dt)$ on $S$, 
given by restriction to $S$ of the canonical contact system on $J^k(\mathbb{R}, G/H)$.
A smooth curve $\gamma: (a, b) \to G/H$ is said to be of \emph{type S\/} 
if $j^k(\gamma)|_t \in S$, for all $t \in (a, b)$. 
\end{definition}

Note that the integral curves of the Pfaffian differential system 
$(\mathcal{I}, dt)$ are the $k$-order jets $j^k(\gamma)$ of curves
$\gamma: (a, b) \to G/H$ that satisfy the differential relation
$j^k(\gamma)|_t \in S$, for all $t \in (a, b)$.

\begin{definition}
A \emph{Frenet system\/} of order $k$ on $G/H$ is a triple $(S, \mathbb{A}, \Phi)$,
where
\begin{itemize}
\item[a)] $S \subset J^k(\mathbb{R}, G/H)$ is a $G$-invariant differential relation of order $k$ 
endowed with the induced contact system $(\mathcal{I}, dt)$,
\item[b)] $\mathbb{A} \in P^h(\mathfrak{g})$,
\item[c)] $\Phi: S \to M$ is a smooth equivariant map from 
$S$ onto an open subset $\Phi(S)$ of $M = G \times \mathbb{A}$, the configuration space,
\end{itemize}
with the properties that
\begin{itemize}
\item If $\gamma: (a, b) \to G/H$ is a smooth curve of type $S$ then
$\Gamma = \Phi \circ j^k(\gamma)$ is an integral curve of the control system $(\mathcal{A}, \omega)$;
\item If $\Gamma: (a, b) \to M$ is an integral curve of $(\mathcal{A}, \omega)$ such that
$\mathrm{Im}(\Gamma) \subset \Phi(S)$, then $\gamma = \pi_{G/H} \circ \Gamma: (a, b) \to G/H$
is a curve of type $S$ and $\Gamma = \Phi \circ j^k(\gamma)$.
\end{itemize}
\end{definition}

The method of moving frames \cite{Griffiths2,Griffiths1} 
gives an algorithmic procedure for the construction of the Frenet systems for curves
of constant type in homogeneous spaces (see \cite{Green,Jensen1,Jensen2}).
We refer the reader to \cite{Griffiths1} for the explicit construction of the Frenet
system of generic curves in the affine space $\mathbb{R}^3$, to \cite{Cartan2,MussoNicolodi}
for the Frenet system of generic curves in $\mathbb{RP}^2$,
to \cite{Musso1,Sulanke} for the Frenet system of generic curves in the conformal $3$-sphere
and to \cite{Musso3} for the Frenet systems of generic Legendrian curves
in the strongly pseudoconvex real hyperquadric $Q^3$ of $\mathbb{CP}^2$.

\begin{definition}
Let $F:= \pi_G \circ \Phi: S \to G$ and 
$K:= \pi_{\mathbb{A}} \circ \Phi: S \to \mathbb{A}$ 
denote the two components of the map $\Phi$.
We call $F$ the \emph{Frenet map\/} and $K$ the \emph{curvature map\/}.

\smallskip

Let $\gamma: (a, b) \to G/H$ be a curve of type $S$,
then $\Gamma = \Phi \circ j^k(\gamma): (a, b) \to M$
is called the \emph{canonical lift\/} of $\gamma$.
The maps
\begin{align*}
F_{\gamma} &: = F \circ j^k(\gamma): (a, b) \to G, \\
K_{\gamma} &: = K \circ j^k(\gamma): (a, b) \to \mathbb{A}
\end{align*}
are called
the \emph{Frenet frame field\/} and the \emph{curvature function\/} of $\gamma$, respectively.
\end{definition}

\begin{definition}
The \emph{generalized arc-length\/} of a curve 
$\gamma: (a, b) \to G/H$ of type $S$ is the smooth function
$s_{\gamma}: (a, b) \to \mathbb{R}$, unique up to a constant, 
such that $ds_{\gamma} = \Gamma^*(\omega)$,
where $\Gamma: (a, b) \to M$ is the canonical lift of $\gamma$. 
Each curve $\gamma \subset G/H$ of type $S$ may be parametrized in such a way that $ds_{\gamma}=dt$.
In this case, we say that the curve $\gamma$ is \emph{normalized\/}.
\end{definition}

\begin{proposition}
Let $(S, \mathbb{A}, \Phi)$ be a Frenet system. 
Then $\Phi(S) = G \times U_{\Phi}$, where $U_{\Phi}$ is an open subset of $\mathbb{A}$.
\label{uphi}
\end{proposition}

\begin{proof}
We set $U_{\Phi} = \pi_{\mathbb{A}}(\Phi(S))$. Thus, $U_{\Phi}$ is an open subset of $\mathbb{A}$ such that
$\Phi (S) \subseteq G \times U_{\Phi}$. Take any $(g_0, Q_0) \in G \times U_{\Phi}$.
Since $Q_0 \in U_{\Phi}$ then there exists $g_1 \in G$ such that $(g_1, Q_0) \in \Phi(S)$.
Now let $\Gamma: ( -\epsilon, \epsilon) \to M$ be an integral curve of the control system 
$(\mathcal{A}, \omega)$ such that $\Gamma(0) = (g_1, Q_0)$.
Since $\Phi(S)$ is an open set we may, 
by restricting the value of $\epsilon$ if necessary, 
assume that $\mathrm{Im}(\Gamma) \subset \Phi(S)$.
Then, the projection of $\Gamma$ onto $G/H$ is a curve
$\gamma: (-\epsilon, \epsilon) \to G/H$ of type $S$ such that $\Gamma = \Phi[j^k(\gamma)]$.
Using the $G$-invariance of $S$ it follows that $g_0 g_1^{-1} \gamma$ is another curve of type $S$.
Thus, from the equivariance of $\Phi$ it follows that $g_0 g_1^{-1} \Gamma(0) = (g_0, Q_0)$
belongs to $\Phi(S)$. This shows that $G\times U_{\Phi} \subseteq \Phi(S)$.
\end{proof}

The elements of the open subset $U_{\Phi}$ may therefore be considered as the 
\lq\lq geometrical inputs\rq\rq\ of the control system $(\mathcal{A}, \omega)$. 
In particular, the curvature function $K$ gives a complete set of local differential invariants for curves of type $S$. 
More precisely, if $\gamma, \widetilde{\gamma}: (a, b) \to G/H$ are normalized curves of type $S$ with 
$K_{\gamma} = K_{\widetilde{\gamma}}$, then $\gamma$ and $\widetilde{\gamma}$ are congruent to one other, 
in the sense that there exists a $g \in G$ such that $g \gamma(t) = \widetilde{\gamma}(t)$, 
for all $t \in (a, b)$. 
Moreover given any smooth map $K: (a, b) \to U_{\Phi} \subset \mathbb{A}$ 
there exists a normalized curve $\gamma: (a, b) \to G/H$ of type $S$, 
unique up to congruence, such that $K_{\gamma} = K$.

If we fix an affine frame $(P, e_1, \dots, e_h)$ of $\mathbb{A}$ and if we let $k^1, \dots, k^h$ be the corresponding coordinates, 
we may identify the configuration space $M$ with $G \times \mathbb{R}^h$.
Thus, we may write $K_{\gamma}=(k^1_{\gamma}, \dots, k^h_{\gamma})$, 
where $k^1_{\gamma}, \dots, k^h_{\gamma}$ are smooth functions that depend on the $k$-jet of $\gamma$.
These functions can be viewed as the \emph{generalized curvatures\/} of $\gamma$.

\subsection{Isotropic curves in $\mathbb{R}^{(2, 1)}$}
\label{sec:R21}
An example that will illustrate our considerations concerns 
variational principles for isotropic curves in $3$-dimensional Minkowski space. 
Let $\mathbb{R}^{(2, 1)}$ denote Minkowski $3$-space endowed with the Lorentzian inner product 
\[
\langle v, w \rangle = - \left( v^1 w^3 + v^3 w^1 \right) + v^2 w^2 =: g_{ij} v^i w^j.
\]
We fix the spatial orientation by requiring that the standard basis $\left( e_1, e_2, e_3 \right)$ is 
positively oriented, and we fix the time orientation defined by the positive light cone
\[
\mathcal{L}^+ = \left\{ v \in \mathbb{R}^{(2, 1)}: \langle v, e_1 + e_3 \rangle < 0 \right\}.
\]
Let $G$ be the restricted Poincar{\'e} group $\mathbb{E}(2, 1)$ 
i.e. the group of isometries of $\mathbb{R}^{(2, 1)}$ that preserve the given orientations.
The group $G$ may conveniently be described as the space of pairs $g = (q, A)$ where 
$q \in \mathbb{R}^{(2, 1)}$ and $A = \left( A_1, A_2, A_3 \right)$ is a $3 \times 3$ matrix 
such that
\[
\det \left( A_1, A_2, A_3 \right) = 1, \qquad 
\langle A_i, A_j \rangle = g_{ij}, \qquad 
A_1, A_3 \in \mathcal{L}^+.
\]
We let $\mathfrak{g}$ denote the Lie algebra of $G$, consisting of all matrices of the form
\[
X(q, v) = 
\begin{pmatrix}
0&0&0&0\\
q^1&v^1_1&v^1_2&0\\
q^2&v^2_1&0&v^1_2\\
q^3&0&v^2_1&-v^1_1
\end{pmatrix}.
\]
We now define the Maurer-Cartan form $\Omega \in \Omega^1(G, \mathfrak{g})$, which takes the form
\[
\Omega = 
\begin{pmatrix}
0&0&0&0\\
\omega^1&\omega^1_1&\omega^1_2&0\\
\omega^2&\omega^2_1&0&\omega^1_2\\
\omega^3&0&\omega^2_1&-\omega^1_1
\end{pmatrix}
\]
such that
\[
dq = \omega^i A_i, \qquad dA_i = \omega^j_i A_j, \qquad i = 1, 2, 3.
\]
Differentiating these relations, we obtain the structure equations
\[
d{\omega}^i = - {\omega}^i_j \wedge {\omega}^j, \qquad 
d{\omega}^i_k = - {\omega}^i_j \wedge {\omega}^j_k, \qquad 
i, k = 1, 2, 3.
\]
Recall that the Maurer-Cartan forms 
$\omega^1, \omega^2, \omega^3, \omega^1_1, \omega^2_1, \omega^1_2$ 
are linearly independent and generate the space, $\mathfrak{g}^*$, 
of left-invariant $1$-forms on $G$.

\begin{definition}
A \emph{null (or isotropic) curve} in {\Mink} is a smooth parametrized curve 
\[
\alpha: \left( a, b \right) \subset \mathbb{R} \rightarrow \Mink
\]
such that $\alpha^{\prime}(t) \in \mathcal{L}^+$ for all $t \in \left( a, b \right)$. 
We shall assume that $\alpha$ is without flex points, in the sense that
\[
\alpha^{\prime}(t) \wedge \alpha^{\prime\prime}(t) \neq 0, \quad \forall t \in \left( a, b \right).
\]
\end{definition}

The linear differential form $\omega_{\alpha} := \Vert \alpha^{\prime\prime}(t) \Vert^{1/2} dt$ is 
nowhere vanishing, and is invariant under changes of parameter and the action of the group $G$. 
Without loss of generality we may assume that $\alpha$ is \emph{normalized}, in the sense that
\[
\Vert \alpha^{\prime\prime}(t) \Vert^{1/2} = 1, \quad \forall t \in \left( a, b \right).
\]
(This condition fixes the parameter $t$ up to an additive constant.) 
The \emph{curvature} of $\alpha$ is defined by
\[
k(t) = - \frac{1}{2} \Vert \alpha^{\prime\prime\prime}(t) \Vert^2, \quad \forall t \in \left( a, b \right).
\]
At each point of the curve we may define the frame $g(t) = \left( \alpha(t), A(t) \right) \in G$ given by
\[
A_1(t) = \alpha^{\prime}(t), \quad
A_2(t) = \alpha^{\prime\prime}(t), \quad
A_3(t) = \alpha^{\prime\prime\prime}(t) + 
\frac{1}{2} \Vert \alpha^{\prime\prime\prime}(t) \Vert^2 \alpha^{\prime}(t).
\]
This frame defines a canonical lift
\[
g: t \in \left( a, b \right) \rightarrow g(t) = \left( \alpha(t), A(t) \right) \in G
\]
of the curve $\alpha$ to the group $G$, referred to as the \emph{Frenet frame field} along $\alpha$. 
An application of the method of moving frames shows that the Frenet frame field is the unique lift 
of $\alpha$ to $G$ with the property that
\[
g^* \left( \Omega \right) = 
\begin{pmatrix}
0&0&0&0\\
1&0&\kappa&0\\
0&1&0&\kappa\\
0&0&1&0
\end{pmatrix} dt.
\]

We illustrate the construction of the Frenet system for isotropic curves in {\Mink}, 
viewed as a homogeneous space of the group $G$. Let 
\[
t, \quad
X = (x^1, x^2, x^3), \quad
X_1 = (x^1_1, x^2_1, x^3_1), \quad
X_2 = (x^1_2, x^2_2, x^3_2), \quad
X_3 = (x^1_3, x^2_3, x^3_3)
\]
be the standard coordinates on the jet space $J^3(\mathbb{R}, \Mink) \cong \mathbb{R} \times \Mink \times \Mink \times \Mink \times \Mink$. 
The differential relation $S \subset J^3(\mathbb{R}, \Mink)$ is defined by
\[
X_1 \in \mathcal{L}^+, \quad
\Vert X_2 \Vert = 1, \quad
\left( X_1, X_2 \right) = \left( X_2, X_3 \right) = 0, \quad
X_1 \wedge X_2 \wedge X_3 \neq 0.
\]
Holonomic sections of $S$ are third order jets $j^3(\alpha)$ of normalized isotropic curves 
$\alpha: \left( a, b \right) \rightarrow \Mink$. We define $\kappa: S \rightarrow \mathbb{R}$ by
\[
\kappa(t, X, X_1, X_2, X_3) = - \frac{1}{2} \Vert X_3 \Vert^2.
\]
$\kappa\left[j^3(\alpha)\right]$ is then the curvature of the isotropic curve $\alpha$. 

The affine space $\mathbb{A} = P + \mathfrak{a} \subset \mathfrak{g}$ is then the straight line
\[
k \in \mathbb{R} \rightarrow Q(k) = e_0 + k e_1 \in \mathfrak{g}
\]
where
\[
e_0 = 
\begin{pmatrix}
0&0&0&0\\
1&0&0&0\\
0&1&0&0\\
0&0&1&0
\end{pmatrix}, \qquad
e_1 = 
\begin{pmatrix}
0&0&0&0\\
0&0&1&0\\
0&0&0&1\\
0&0&0&0
\end{pmatrix}.
\]
It is convenient to identify the configuration space $M = G \times \mathbb{A}$ with $G \times \mathbb{R}$ by 
means of the map 
\[
\left( g, e_0 + k e_1 \right) \in G \times \mathbb{A} \rightarrow \left( g, k \right) \in G \times \mathbb{R}.
\]
With this identification at hand, the linear control system $(\mathcal{A}, \omega)$ is generated by the linear differential forms
\[
\eta^1 = {\omega}^1_2 - k \omega, \quad
\eta^2 = {\omega}^1_1, \quad
\eta^3 = {\omega}^2_1 - \omega, \quad
\eta^4 = {\omega}^2, \quad
\eta^5 = {\omega}^3.
\]
along with independence condition
\[
\omega = {\omega}^1.
\]

We now consider a smooth curve $\Gamma: \left( a, b \right) \rightarrow M$ and let
\[
g: t \in \left( a, b \right) \rightarrow \left( \alpha(t), A(t) \right) \in G, \qquad
k: t \in \left( a, b \right) \rightarrow k(t) \in \mathbb{A},
\]
be the two components of $\Gamma$. Then $\Gamma$ is an integral curve of the control system $(\mathcal{A}, \omega)$ 
if and only if $g: \left( a, b \right) \rightarrow G$ is the Frenet field along the isotropic curve 
$\alpha: \left( a, b \right) \rightarrow \Mink$, and $k$ is the curvature of the curve $\alpha$. 
The mapping $\Phi: S \rightarrow M = G \times \mathbb{R}$ linearizing the differential relation $S$ is defined by
\[
S: (t, X, X_1, X_2, X_3) \in S \rightarrow 
\left( \left( X; X_1, X_2, X_3 + \frac{1}{2} \Vert X_3 \Vert^2 X_1 \right), - \frac{1}{2} \Vert X_3 \Vert^2 \right) \in M.
\]

\begin{remark}
Using the structure equations for $\mathbb{E}(2, 1)$ we find that
\begin{subequations}
\begin{align}
d\omega &= \left( \kappa \eta^4 - \eta^2 \right) \wedge \omega - \eta^1 \wedge \eta^4,\\
d\eta^1 &= - \pi \wedge \omega + \eta^1 \wedge \eta^2 + \kappa \eta^1 \wedge \eta^4,\\
d\eta^2 &= \left( \kappa \eta^3 - \eta^1 \right) \wedge \omega - \eta^1 \wedge \eta^3,\\
d\eta^3 &= \left( 2 \eta^2 - \kappa \eta^4 \right) \wedge \omega + \eta^1 \wedge \eta^4 + \eta^2 \wedge \eta^3,\\
d\eta^4 &= \left( \kappa \eta^5 - \eta^4 \right) \wedge \omega - \eta^1 \wedge \eta^5,\\
d\eta^5 &= \eta^4 \wedge \omega  + \eta^2 \wedge \eta^5 - \eta^3 \wedge \eta^4,
\end{align}
\label{structure}
\end{subequations}
where
\[
\pi = d\kappa + \kappa^2 \eta^4.
\]
\end{remark}

\subsection{Coadjoint action of $\mathbb{E}(2, 1)$}
\label{sec:coadj}
For later convenience, we now discuss the coadjoint action of $\mathbb{E}(2, 1)$ on $\mathfrak{e}(2, 1)^*$,
the dual of its Lie algebra. Our discussion follows the discussion of the coadjoint representation of $\mathbb{E}(3)$
given in Guillemin and Sternberg \cite{GuilleminSternberg1}.

Using the Maurer-Cartan forms, we identify $\mathfrak{g}^*$ with $\Mink \oplus \Mink$ by means of the map
\[
(p, v) \in \Mink \oplus \Mink \rightarrow 
p_i \omega^i - v^1 \omega^2_1 + v_2 \omega^1_1 + v_3 \omega^1_2 \in \mathfrak{g}^*.
\]
The coadjoint action of $G$ on $\mathfrak{g}^*$ then takes the form
\begin{equation}
g \cdot \left( p, v \right) = \left( A p, A v - \left( A p \right) \times q \right),
\label{coadjointaction}
\end{equation}
for all $g = \begin{pmatrix}1&0\\q&A\end{pmatrix} \in G = \mathbb{E}(2, 1)$,
where $\times$ denotes the vector cross product
\[
\langle v \times w, u \rangle = \det (v, w, u ), \qquad \forall v, w, u \in \Mink.
\]
We now define the map
\[
C: \left( p, v \right) \in \mathfrak{g}^* \rightarrow \left( \Vert p \Vert^2, \langle p, v \rangle \right) \in \mathbb{R}^2,
\]
the components which, $C_1$ and $C_2$, generate the space of Casimir functions. We recall the following standard material:
\begin{itemize}
\item
Let $G$ be a Lie group, and $\mathfrak{g}^*$ the dual of the Lie algebra of $G$.
Let $\mu \in \mathfrak{g}^*$. The \emph{isotropy group} of $G$ at $\mu$ is the
closed subgroup of $G$ defined by
\[
G_{\mu} := \left\{ g \in G: Ad^*(g) \mu = \mu \right\} =
\left\{ g \in G: \langle \mu; Ad(g^{-1}) A \rangle = \langle \mu; A \rangle, \forall A \in \mathfrak{g} \right\}.
\]
\item
The Lie algebra of $G_{\mu}$ is
\[
\mathfrak{g}_{\mu} = \left\{ A \in \mathfrak{g}: ad^*(A) \mu = 0 \right\} =
\left\{ A \in \mathfrak{g}: \langle \mu; [A, B] \rangle =0, \forall B \in \mathfrak{g} \right\}
\]
\item
The \emph{rank\/} of the group $G$ is defined as
\[
\rm{rank}(G) = \rm{inf} \{ \rm{dim}(\mathfrak{g}_{\mu}): \mu \in \mathfrak{g}^* \}.
\]
\item
An element $\mu \in \mathfrak{g}^*$ is \emph{regular\/} if
$\rm{dim}(\mathfrak{g}_{\mu}) = \mathrm{rank}(G)$,
otherwise $\mu$ is a \emph{singular\/} element of $\mathfrak{g}^*$.
The set of regular elements of $\mathfrak{g}^*$ will be denoted by $\mathfrak{g}^*_r$,
while $\mathfrak{g}^*_s$ will denote the set of singular elements.
\item
By a theorem of Dixmier (cf. \cite{Dixmier,FomenkoTrofimov}),
the isotropy group $G_{\mu}$ and the isotropy Lie algebra $\mathfrak{g}_{\mu}$
of a regular element $\mu \in \mathfrak{g}^*_r$ are Abelian.
\end{itemize}

\medskip

In the case of $\mathbb{E}(2, 1)$, $\mathfrak{g}^*_r$ is the open subset of $\mathfrak{g}^*$ consisting of elements
\[
\mathfrak{g}^*_r = \left\{ \left( p, v \right) \in \mathfrak{g}: p \neq 0 \right\}.
\]
The co-adjoint orbit, $\mathcal{O}_{(p_0, v_0)}$, through a regular element $(p_0, v_0) \in \mathfrak{g}^*_r$ is therefore the
four-dimensional sub-manifold
\[
\mathcal{O}_{(p_0, v_0)} = \left\{ \left( p, v \right) \in \Mink \oplus \Mink:
\Vert p \Vert^2 = \Vert p_0 \Vert^2,
\langle p_0, v_0 \rangle = \langle p, v \rangle \right\}.
\]
There are three types of regular orbit:
\begin{itemize}
\item orbits of positive type: $\mathcal{O}_{(p_0, v_0)}$ with $C_1 = \Vert p_0 \Vert^2 > 0$;
\item orbits of negative type: 
$\mathcal{O}_{(p_0, v_0)}$ with $C_1 = \Vert p_0 \Vert^2 < 0$;
\item orbits of null type: $\mathcal{O}_{(p_0, v_0)}$ with $C_1 = \Vert p_0 \Vert^2 = 0$;
\end{itemize}
The orbits of negative and null type also break into sub-classes according to whether 
$p_0$ is future-directed, with $\langle p_0, e_1 + e_3 \rangle < 0$, or 
past-directed, with $\langle p_0, e_1 + e_3 \rangle > 0$.

\section{Variational problems}
\label{sec:variational}
\subsection{Non-degenerate invariant variational problems}

\begin{definition}
Given an affine subspace $\mathbb{A} \in P^h (\mathfrak{g})$, 
an \emph{invariant Lagrangian of type $\mathbb{A}$\/}
is a smooth real-valued function $L: \mathbb{A} \to \mathbb{R}$.
\end{definition}

An invariant Lagrangian $L$ gives
rise to a variational problem defined on the integral curves of
the linear control system $(\mathcal{A}, \omega)$. From this point of view
the Lagrangian $L$ is considered as a \emph{cost function\/}. It
is then an optimal control problem to minimize the cost
\[
\mathcal{L}: \Gamma \to \int _{\gamma}\Gamma^*(L \, \omega)
\]
among the integral curves of $(\mathcal{A}, \omega)$. 
If $(\mathcal{A}, \omega)$ comes from a Frenet system 
$(S, \mathbb{A}, \Phi)$ on the homogeneous space $G/H$ 
then the Lagrangian $L$ defines a geometric action functional 
$\widetilde{\mathcal{L}}: \mathcal{S} \to \mathbb{R}$ 
acting on the space of the normalized curves of type $S$:
\[
\widetilde{\mathcal{L}}: 
\gamma \in \mathcal{S} \to \int_{\gamma} L \left( K[j^k(\gamma)(t)] \right) \, dt.
\]
Note that the geometric action functional $\widetilde{\mathcal{L}}$ 
depends only on the generalized curvatures of $\gamma$.

\begin{example}
The simplest invariant variational problem for a
Frenet systems is the \emph{arc-length functional\/}, 
which is defined by a constant Lagrangian 
(see \cite{Cartan2,Musso1,Musso3,MussoNicolodi,Griffiths1} for more
details about the arc-length functionals for generic curves in the
conformal and pseudoconformal three-dimensional sphere, in the
real projective plane and in the affine plane). 
Another typical example of an invariant Lagrangian is the 
\emph{Kirchhoff variational problem\/} 
for the Frenet system of generic curves in $\mathbb{R}^3$, 
defined by the action functional
\[
\mathcal{L}: \gamma \subset \mathbb{R}^3 \to 
\int_{\gamma} \left( \kappa^2(u) + a \tau(u) \right) du.
\] 
The extremal curves are the canonical lifts of the Kirchhoff elastic rods of $\mathbb{R}^3$. 
When $a=0$, we get the total squared curvature functional. 
Other examples of invariant variational problems for curves in $\mathbb{R}^3$ 
have been considered in ref.~\cite{LangerSinger}.
\end{example}

Given an invariant Lagrangian $L: \mathbb{A} \to \mathbb{R}$,
we construct the corresponding affine sub-bundle ${\widetilde Z} \subset T^*(M)$
over the configuration space $M = G \times \mathbb{A}$.
The fibre of ${\widetilde Z}$ over the point $(g, Q) \in M$ is given by the affine space
\[
\left. {\widetilde Z} \right|_{(g, Q)} = \{ \eta \in \mathfrak{g}^*: \langle \eta; Q \rangle = L(Q) \}.
\]
Note that ${\widetilde Z}$ is of the form $G \times \mathcal{{\widetilde Z}}$, where
\[
\mathcal{{\widetilde Z}} = \{ (Q, \eta) \in \mathbb{A} \times \mathfrak{g}^*: \langle \eta; Q \rangle = L(Q) \}.
\]
The Liouville one-form $\psi$ is given by
\[
\psi|_{(g, Q, \eta)} = \pi^*(\eta)|_{(g, Q, \eta)}, \quad \forall (g, Q, \eta) \in {\widetilde Z},
\]
where $\pi: G \times \mathcal{{\widetilde Z}} \to G$ denotes the projection onto the first factor.

\begin{remark}
Pick a basis $(e_0, e_1, \dots, e_h, e_{h+1}, \dots, e_n)$ of $\mathfrak{g}$
such that
\[
P = e_0, \quad \mathfrak{a} = \mathrm{span}(e_1, \dots, e_h), \quad \omega = \theta^0,
\]
where $(\theta^0, \dots, \theta^n)$ is the dual basis of $\mathfrak{g}^*$. 
We use the following index range: 
$i, j = 1, \dots, h$, $a, b = h+1, \dots, n$. 
The map
\[
(g, k, \lambda) \in G \times \mathbb{R}^h \times \mathbb{R}^n 
\to (g, e_0 + k^j e_j, L(e_0 + k^j e_j) \omega + 
\lambda_j (\theta^j - k^j \omega) + 
\lambda_a \theta^a) \in {\widetilde Z}
\]
gives an explicit identification between $G \times \mathbb{R}^h \times \mathbb{R}^n$ and ${\widetilde Z}$. 
With this identification at hand, the tautological $1$-form can be written as:
\[
\psi = (L(k^1, \dots, k^h) - k^j \lambda_j) \omega + \lambda_j \theta^j + \lambda_a \theta^a.
\]
\end{remark}

\begin{definition}
An invariant Lagrangian $L: \mathbb{A} \to \mathbb{R}$ is said to
be \emph{regular\/} if the corresponding variational problem $(\mathcal{A}, \omega, L)$ is regular 
i.e. if the Cartan system of $\Psi = d\psi$, with the independence condition $\omega$, 
is reducible (see Definition~\ref{regular}). 
For a regular Lagrangian we denote by $Y \subset {\widetilde Z}$ the momentum space of the variational problem
$(\mathcal{A}, \omega, L)$.
\end{definition}

\begin{remark}
We have seen that all the derived systems of $(\mathcal{A}, \omega)$
have constant rank. This implies that the extremal curves of a
regular invariant variational problem are the projections of the
integral curves of the Euler-Lagrange system on $Y$ (c.f. \cite{Bryant}).
\end{remark}

\begin{proposition}
Let $L: \mathbb{A} \to \mathbb{R}$ be a regular Lagrangian with momentum space $Y$. 
Then $Y = G \times \mathcal{F}$,
where $\mathcal{F}$ is an immersed submanifold of $\mathbb{A} \times \mathfrak{g}^*$.
\end{proposition}

\begin{proof}
First we claim that the momentum space, $Y$, is $G$-invariant.
To show this, for any $g \in G$, 
we consider the submanifold $g \cdot Y \subset {\widetilde Z}$. 
The $G$-invariance of the exterior differential forms $\psi$, $\Psi$ and $\omega$ 
implies that left translation $L_g: {\widetilde Z} \to {\widetilde Z}$ sends integral elements of
$(\mathcal{C}(\Psi), \omega)$ into integral elements of $(\mathcal{C}(\Psi), \omega)$.
Hence, for every point $p \in g \cdot Y$, there exists an integral
element of $(\mathcal{C}(\Psi), \omega)$ tangent to $g \cdot Y$. 
Since the momentum space $Y$ is maximal with respect to this property, 
it follows that $g \cdot Y \subseteq Y$. Thus the group $G$ acts on $Y$. 
Since this action is free and proper, the quotient space $\mathcal{F} = Y/G$
exists as a manifold. The natural projection $\pi: Y \to \mathcal{F}$ is
constant along the fibers of the map $(g, Q, \eta) \in Y \to
(Q, \eta) \in \mathbb{A} \times \mathfrak{g}^*$. Thus it induces a smooth
one-to-one immersion $j: \mathcal{F} \to \mathbb{A} \times \mathfrak{g}^*$. We conclude the
proof by observing that the map $(\mathrm{id}, j): G \times \mathcal{F} \to Y$
is a smooth diffeomorphism.
\end{proof}

\begin{definition}
We call $\mathcal{F}$ the \emph{phase space\/} of the system. 
Note that a point $p \in \mathcal{F}$ is of the form 
$p = (Q, \eta)$, where $Q \in \mathbb{A}, \eta \in \mathfrak{g}^*$. 
We define the maps 
\begin{align*}
\Lambda&: (Q, \eta) \in \mathcal{F} \to \eta \in \mathfrak{g}^*,\\
\mathcal{H}&: (Q, \eta) \in \mathcal{F} \to Q \in \mathbb{A} \subset \mathfrak{g}.
\end{align*}
We refer to $\Lambda$ as the \emph{Legendre transform\/} and $\mathcal{H}$ as the \emph{Hamiltonian\/}.
Let $F(p) := T_p(\mathcal{F}) \subset \mathfrak{a} \oplus \mathfrak{g}^*$
be the tangent space of $\mathcal{F}$ at $p$. We then define
\[
R(p) := d\Lambda|_p [F(p)] \subset \mathfrak{g}^*, \quad
S(p) := d\mathcal{H}|_p [F(p)] \subset \mathfrak{a}, \quad \forall p \in \mathcal{F}.
\]
\end{definition}

\begin{definition}
A regular invariant Lagrangian $L: \mathbb{A} \to \mathbb{R}$ is said to be \emph{non-degenerate\/}
if the momentum space $Y$ is odd-dimensional, of dimension $2m+1$,
and if the restriction of the canonical $2$-form $\Psi$ to $Y$, $\Psi_Y$, 
has the property that $\omega \wedge \left( \Psi_Y \right)^m$ is non-vanishing.
\end{definition}

Examples of invariant non-degenerate variational problems include the
total squared curvature functional in two and three-dimensional
space forms \cite{Griffiths1,BryantGriffiths}, the Kirchhoff
variational problem in $\mathbb{R}^3$, the Poincar\'e and the Delaunay
functionals \cite{Griffiths1,Hsu,Musso2}, the projective, the
conformal and the pseudoconformal arc-length functionals (cf.
\cite{Cartan2,Musso1,Musso3}).

\medskip

Given a non-degenerate variational problem,
it follows that $\omega \wedge \left( \Psi_Y \right)^m$ defines a volume form on $Y$,
and that $\Psi_Y$ is of maximal rank on $Y$.
Therefore there exists a unique vector field $\xi \in \mathfrak{X}(Y)$ such that
$i_{\xi} \left( \Psi_Y \right) = 0$ and $\omega(\xi) =1$.

\begin{definition}
$\xi$ is the \emph{characteristic vector field\/} of the non-degenerate variational problem $(\mathcal{A}, \omega, L)$.
\end{definition}

If $(\mathcal{A}, \omega, L)$ is non-degenerate then the Euler-Lagrange system  is simply the Cartan system of the
canonical $2$-form restricted to the momentum space: $\mathcal{E} = \mathcal{C}( \Psi_Y)$.
Therefore, for such variational problems,
the integral curves of the Euler-Lagrange system
are the integral curves of the characteristic vector field $\xi$
(see ref.~\cite{Griffiths1}). We therefore have

\begin{thm}
Let $\Gamma: (a, b) \to Y$ be an integral curve of the characteristic vector field $\xi$
of a non-degenerate variational problem $(\mathcal{A}, \omega, L)$.
Then $\gamma = \pi_M \circ \Gamma: (a, b) \to M$ is a critical point of the action functional $\mathcal{L}$.
\end{thm}
\vskip .2cm

\begin{proposition}
If $L$ is non-degenerate then the Legendre transform 
$\Lambda: \mathcal{F} \to \mathfrak{g}^*$ is an immersion.
\end{proposition}

\begin{proof}
Let $\xi$ be the characteristic vector field of the momentum space.
The Liouville form $\psi$, the canonical two-form $\Psi$ 
and the independence condition $\omega$ are $G$-invariants, 
therefore the characteristic vector field is also $G$-invariant. 
Since $\xi|_{(g, p)} \in T_{(g, p)}Y \cong  \mathfrak{g} \oplus F(p) \subset \mathfrak{g} \oplus \mathfrak{a} \oplus \mathfrak{g}^*$,
this implies that there exist smooth maps
$A_{\xi}: \mathcal{F} \to \mathfrak{g}$ and
$\Phi_{\xi}: \mathcal{F} \to \mathfrak{a} \oplus \mathfrak{g}^*$ with the property that
\[
\xi|_{(g, p)} = A_{\xi}(p)|_{g} + \Phi_{\xi}(p), \quad \forall (g, p) \in Y,
\]
where $\Phi_{\xi}(p) \in F(p)$, for all $p \in \mathcal{F}$.
Since $\xi$ satisfies the transversality condition
$1 = \omega(\xi) = \langle \omega; A_{\xi} \rangle$, 
then $A_{\xi}: \mathcal{F} \to \mathfrak{g}$ is a nowhere vanishing function. 
If we now consider $\{ 0 \} \oplus \ker [ d\Lambda|_p \,] \subset \mathfrak{g} \oplus F(p)$ 
then it is simple to check that every such vector lies in the kernel of the canonical $2$-form $\Psi$. 
Since this null-distribution is generated by $\xi$, we therefore have 
\[
\{ 0 \} \oplus \ker [ d\Lambda|_p \,] \subset \mathrm{span}[A_{\xi}(p) + \Phi_{\xi}(g)].
\]
Since $A_{\xi}$ is non-vanishing, however, 
this holds if and only if $\ker \left[ d\Lambda|_p \, \right] = \{ 0 \}$.
\end{proof}
\vskip .2cm

From now on we will assume that the Legendre transform $\Lambda$ 
is a one-to-one immersion, so that the phase space, $\mathcal{F}$, 
can be considered as a submanifold (not necessarily embedded) of $\mathfrak{g}^*$. 
Consequently, we will think of the momentum space 
as an immersed submanifold of $G \times \mathfrak{g}^*$. 
The notation introduced in the preceding paragraphs can then be simplified as follows:
\medskip
\begin{itemize}
\item the Legendre map $\Lambda$ is the inclusion of $\mathcal{F}$ into $\mathfrak{g}^*$;
\item the tangent space $F(\eta)$ of $\mathcal{F}$ at
$\eta \in \mathcal{F}$ is a linear subspace of $\mathfrak{g}^*$ and $R(\eta) = F(\eta)$;
\item the tangent space $T_{(g, \eta)}(Y)$ is identified with
$\mathfrak{g} \oplus F(\eta) \subset \mathfrak{g} \oplus \mathfrak{g}^*$;
\item the Liouville form and the canonical two-form on $Y$ are the
restrictions to $Y$ of the Liouville form and the standard
symplectic form on $T^*(G)$;
\item the characteristic vector field $\xi$ can be written as
\[
\xi|_{(g, \eta)} = \left. A_{\xi}(\eta) \right|_g + \Phi_{\xi}(\eta), \quad
\forall (g, \eta) \in Y,
\]
where $A_{\xi}: \mathcal{F} \to \mathfrak{g}$ and
$\Phi_{\xi}: \mathcal{F} \to \mathfrak{g}^*$ 
are smooth functions such that 
$\Phi_{\xi}(\eta) \in F(\eta)$, for all $\eta \in \mathcal{F}$.
\end{itemize}
From now on, we will adhere to these simplifications.

With this notation at hand, we may use the left-invariant trivialization of $T(G)$ to identify the tangent space
\[
T_{(g, \eta)}(Y) \cong 
T_g G \oplus T_{\eta}\mathcal{F} \cong 
\mathfrak{g} \oplus F(\eta) \subset 
\mathfrak{g} \oplus \mathfrak{g}^*.
\]
We then have the explicit isomorphism
\[
A + v \in \mathfrak{g} \oplus F(\eta) \rightarrow  \left. A \right|_g + v \in T_{(g, \eta)}(Y),
\]
where $A \in \mathfrak{g} = T_{id}(G)$ and $\left. A \right|_g = \left( L_g \right)_* A \in T_g (G)$. 
With this identification, 
the Liouville form $\psi$ becomes the cross-section of $T^*(Y)$ defined by
\begin{equation}
\left. \psi \right|_{(g, \eta)} (A + v) = 
\langle \eta, A \rangle, \quad 
\forall (g, \eta) \in Y, \quad
\forall A + v \in \mathfrak{g} \oplus F(\eta).
\label{Liouvilleform}
\end{equation}
Then, from the standard formula
\[
d\psi(X, Y) = 
\frac{1}{2} \left\{ X[\psi(Y)] - Y[\psi(X)] - \psi([X, Y]) \right\}, 
\]
it follows that the canonical $2$-form $\Psi = d\psi \in \Omega^2(Y)$ takes the form
\begin{equation}
\left.\Psi\right|_{(g, \eta)} (A+v; B+w) = 
- \frac{1}{2} \langle w; A \rangle 
+ \frac{1}{2}\langle ad^*(A) \eta + v; B \rangle, 
\label{canonicaltwoform}
\end{equation}
for all $\eta \in \mathcal{F}$ and for all 
$A + v, \, B + w \in \mathfrak{g} \oplus F(\eta)$.

\begin{definition}
Given a left-invariant $1$-form $\mu \in \mathfrak{g}^*$, 
let $\mathcal{O}(\mu) \subset \mathfrak{g}^*$ be the coadjoint orbit passing through $\mu$, 
and let $O(\mu):= ad^*(\mathfrak{g}) \mu \subset \mathfrak{g}^*$ 
denote the tangent space to the orbit $\mathcal{O}(\mu)$ at $\mu$. 
The \emph{linearized phase portrait\/} of the point $\eta \in \mathcal{F}$ 
is the linear subspace $\Pi(\eta): = F(\eta) \cap O(\eta)$ of $\mathfrak{g}^*$. 
The subset $\mathcal{P}(\mu) = \mathcal{F} \cap \mathcal{O}(\mu)$ 
is referred to as the \emph{phase portrait\/} of $\mu \in \mathfrak{g}^*$. 
\end{definition}

The following result shows that the characteristic vector field $\xi$ 
may be written in terms of the Hamiltonian $\mathcal{H}$:
\begin{thm}
The characteristic vector field $\xi$ is given by 
\begin{equation}
\xi|_{(g, \eta)} = \left. \mathcal{H}(\eta) \right|_g - ad^*[\mathcal{H}(\eta)] \eta, \quad 
\forall (g, \eta) \in Y. 
\label{characteristic}
\end{equation}
\end{thm}

\begin{proof}
Given a point $\eta \in \mathcal{F}$, we set
\[
\mathrm{Ann}\!\left( F(\eta) \right) = \{ A \in \mathfrak{g}: \langle v; A \rangle = 0, \forall v \in F(\eta) \}
\]
and let
\[
\rho(\eta): \mathrm{Ann}\!\left( F(\eta) \right) \to O(\eta)
\]
be the linear map
\begin{equation}
\rho(\eta): A \in \mathrm{Ann}\!\left( F(\eta) \right) \to ad^*(A) \eta \in O(\eta).
\label{rhodef}
\end{equation}
It then follows from equations~\eqref{Liouvilleform} and~\eqref{canonicaltwoform} that 
a tangent vector $A + v \in \mathfrak{g} \oplus F(\eta)$ to the
momentum space $Y$ at the point $(g, \eta)$ belongs to the kernel of
$\Psi$ if and only if
\begin{equation}
A \in \rho(\eta)^{-1} \left( \Pi(\eta) \right), \quad
v = - \rho(\eta) A. 
\label{kerPsi}
\end{equation}

We now let $(g_0, \eta_0) \in Y$ and let $\Gamma: (-\epsilon, \epsilon) \to Y$ 
be the integral curve of the characteristic vector field $\xi$ 
with initial condition $\Gamma(0)=(g_0, \eta_0)$. 
We write $\Gamma(t)=(g(t), \eta(t))$, 
where $g: (-\epsilon, \epsilon) \to G$ and 
$\eta: (-\epsilon, \epsilon) \to \mathcal{F}$ 
are smooth maps such that
\[
g(t)^{-1} g'(t) dt = g^*(\Theta)|_t, \quad
g^{-1}(t)g'(t) = A_{\xi}[\eta(t)], \quad 
g(0)=g_0, \quad
\eta(0)=\eta_0.
\]
On the other hand 
\footnote{It is a general fact that if $\pi: Y \to M$ 
is the momentum space of a regular variational problem 
$(\mathcal{I}, \omega, L)$ on the configuration space $M$ 
and if $\Gamma: (a, b) \to Y$ is an integral curve of the Euler-Lagrange system, 
then $\gamma = \pi \circ \Gamma$ is an integral curve of the 
Pfaffian differential system $(\mathcal{I}, \omega)$ on $M$ (see \cite{Griffiths1}).}
\[
t \in (-\epsilon, \epsilon) \to (g(t), \mathcal{H}[\eta(t)]) \in G \times \mathbb{A} = M
\]
is an integral curve of the linear control system $(\mathcal{A}, \omega)$. 
We then have
\[
\left. g^*(\Theta) \right|_t = 
\mathcal{H}[\eta(t)] \left. g^*(\omega) \right|_t = 
\mathcal{H}[\eta(t)] \left. dt \right|_t.
\]
Therefore, we conclude that
\[
A_{\xi}[\eta(t)] = \mathcal{H}[\eta(t)], \quad \forall t \in (-\epsilon, \epsilon). 
\]
Since $\xi$ belongs to the kernel of $\Psi$, we conclude from equation~\eqref{kerPsi} that
\[
\Phi_{\xi}[\eta(t)] = - ad^*[\mathcal{H}(\eta(t))] \eta(t), \quad 
\forall t \in (-\epsilon, \epsilon).
\]
This yields the required result.
\end{proof}

\begin{definition}
The \emph{phase flow\/} is the flow of the vector field 
$\Phi_{\xi}: \eta \in \mathcal{F} \to - ad^*[\mathcal{H}(\eta)] \eta \in \mathfrak{g}^*$.
\end{definition}

\begin{remark}
We use the notation $\phi_{\xi}: \mathcal{D} \subset \mathbb{R} \times \mathcal{F} \to \mathcal{F}$
to indicate the phase flow. We observe the following facts:
\begin{itemize}
\item the domain of definition $\mathcal{D}$ of the phase flow is of the form
\[
\mathcal{D} = \{ (t, \eta) \in \mathbb{R} \times \mathcal{F}: t \in (\epsilon^{-}(\eta), \epsilon^{+}(\eta)) \}
\]
where $\epsilon^{-}: \mathcal{F} \to \mathbb{R}^{-} \cup \{-\infty\}$ and
$\epsilon^{+}: \mathcal{F} \to \mathbb{R}^{+} \cup \infty$;
\item For every $\eta \in \mathcal{F}$, the curve
\[
\phi_{\eta}: (\epsilon^{-}(\eta), \epsilon^{+}(\eta)) \to \phi_{\xi} (t, \eta) \in \mathcal{F}
\]
is the maximal integral curve of $\Phi_{\xi}$ with the initial condition $\phi_{\eta}(0) = \eta$;
\item $\Phi_{\xi}(\eta) \in \Pi(\eta)$ and $\phi_{\eta}(t) \in \mathcal{P}(\eta)$,
for every $\eta \in \mathcal{F}$ and every $t \in (\epsilon^{-}(\eta), \epsilon^{+}(\eta))$;
\item if we fix a point $(g_0, \eta_0) \in Y = G \times \mathcal{F}$,
then the maximal integral curve of the characteristic vector field $\xi$ with the initial condition $(g_0, \eta_0)$ is given by
\[
\Gamma_{(g_0, \eta_0)}: t \in (\epsilon^{-}(\eta_0), \epsilon^{+}(\eta_0)) \to
(h_{(g_{0}, \eta_{0})}(t), \phi_{\eta_0}(t)) \in G \times \mathcal{F},
\]
where $h_{(g_0, \eta_0)}$ is the (unique) solution of the equation
\[
h^{-1} h' = \mathcal{H}[\phi_{\eta_0}(t)], \quad h(0) = g_0;
\]
\item We set
$\widetilde{\mathcal{D}} = \{ (t; (g, \eta)) \in \mathbb{R} \times Y: t \in (\epsilon^{-}(\eta), \epsilon^{+}(\eta)) \}$.
The flow $\Gamma$ of the characteristic vector field $\xi$ is the local $1$-parameter group of transformations
$\Gamma: \widetilde{\mathcal{D}} \subset \mathbb{R} \times Y \to Y$
given by
\[
\Gamma(t, g, \eta) = (h_{(g, \eta)}(t), \phi(t, \eta)), \quad \forall (t; (g, \eta)) \in \widetilde{\mathcal{D}}.
\]
\end{itemize}
\end{remark}

\begin{remark}
The phase flow $\phi_{\xi}: \mathcal{D} \to \mathcal{F}$ satisfies the 
\emph{Euler equation\/}
\begin{equation}
\left. \frac{\partial \phi_{\xi}}{\partial t} \right|_{(t, \eta)} = - 
ad^*\left[\mathcal{H}[\phi_{\xi}(t, \eta)]\right] \phi_{\xi}(t, \eta), \quad
\phi_{\xi}(0, \eta) = \eta, \quad \forall (t, \eta) \in \mathcal{D}.
\label{euler}
\end{equation}
If there exists a $G$-equivariant isomorphism $\mathfrak{g} \cong \mathfrak{g}^*$, 
then we can identify $\mathfrak{g}$ and $\mathfrak{g}^*$. 
(For example, if $G$ is semisimple, then take the pairing defined by
the Killing form of $\mathfrak{g}$.) 
Using this identification, 
the Euler equation can be written in the Lax form
\[
\left. \frac{\partial \phi_{\xi}}{\partial t} \right|_{(t, \eta)} = - 
\left[ \vphantom{|^|} \mathcal{H}[\phi_{\xi}(t, \eta)], \phi_{\xi}(t, \eta) \right], \quad
\phi_{\xi}(0, \eta) = \eta, \quad \forall (t, \eta) \in \mathcal{D}.
\]
\end{remark}

\begin{definition}
We denote by $\mathcal{F}_s = \{\eta \in \mathcal{F}: \Phi_{\xi}(\eta) = 0 \}$ 
the set of all fixed points of the phase flow and by
$\mathcal{F}_r$ the complement of $\mathcal{F}_s$.
We call $\mathcal{F}_s$ and $\mathcal{F}_r$ the 
\emph{singular\/} and the \emph{regular\/} parts of the phase space, respectively. 
We call $\Sigma = G \times \mathcal{F}_s \subset Y$ the \emph{bifurcation set\/} 
and refer to $Y_r = Y \setminus \Sigma$ as the 
\emph{regular part\/} of the momentum space. 
The intersection $\mathcal{P}_r(\mu) = \mathcal{F}_r \cap \mathcal{O}(\mu) \subset \mathcal{P}(\mu)$ 
is called the \emph{regular part\/} of the phase portrait $\mathcal{P}(\mu)$. 
The connected component $\widetilde{\mathcal{P}}(\mu)$ of $\mathcal{P}_r(\mu)$ containing $\mu$ 
is referred to as the \emph{reduced phase portrait\/} of $\mu$.
\end{definition}

The following result, which may be verified by applying the uniqueness theorems for ordinary differential equations,
characterises integral curves of the characteristic vector field that intersect the bifurcation set:
\begin{proposition}
Let $p = (g, \eta) \in \Sigma$ be a point of the bifurcation set.
The integral curve $\Gamma_{\xi}(-, p): \mathbb{R} \to Y$
of the characteristic vector field $\xi$ passing through $p$
is the orbit of the one-parameter subgroup generated by $\mathcal{H}(\eta)$:
\[
\Gamma_{\xi}(t, p) =
\left( \mathrm{Exp} \left( \mathcal{H}(\eta) t \right) g, \eta \right), 
\quad \forall t \in \mathbb{R}.
\]
\end{proposition}

This result implies that if $(\mathcal{A}, \omega)$ comes from a Frenet system of
curves in $G/H$, then the curve $\gamma \subset G/H$ of type $S$
that corresponds to $\Gamma_{\xi}(-, p)$, where $p \in \Sigma$, has
constant curvature (i.e. $K_{\gamma} = \mathrm{constant})$.

Since this result completely characterises the behaviour of integral curves that intersect the bifurcation set $\Sigma$,
we shall henceforth restrict our attention to the regular parts of the phase space and momentum space.
Therefore, to simplify the notation, $\mathcal{F}$, $Y$ and $\mathcal{P}(\mu)$ will be used to denote the regular parts
of the phase space, the momentum space and the phase portraits, respectively.

\subsection{The Poincar{\'e} variational problem for isotropic curves in {\Mink}}
We now return to our example of isotropic curves in {\Mink} considered in Section~\ref{sec:R21}. 
Let $m$ be a non-zero constant and consider the variational problem on the space $\mathcal{V}$ of 
parametrized integral curves 
$\Gamma: t \in \left( a, b \right) \rightarrow \left( g(t), k(t) \right) \in G \times \mathbb{R}$ 
of the Pfaffian system $(\mathcal{A}, \omega)$ defined by the action functional
\[
\mathcal{L}_m: \Gamma \in \mathcal{V} \rightarrow 
\int_{\Gamma} \left( 1 + m k \right) \omega.
\]
Geometrically, this amounts to an analogue of the Poincar{\'e} variational problem where we minimize 
the arc-length functional (defined by the integral of the canonical line-element of the null curve) 
amongst normalized null curves $\alpha \subset \Mink$ subject to the additional constraint that the 
integral of the curvature $k$ along the curve be held constant. 

The affine sub-bundle ${\widetilde Z} \subset T^*(M)$ is given by $M \times \mathcal{{\widetilde Z}}$ where 
$\mathcal{\widetilde Z} \subset \mathfrak{g} \oplus \mathfrak{g}^*$ is the submanifold consisting of all 
$(Q(k), \eta) \in \mathbb{A} \oplus \mathfrak{g}^*$ such that $\langle \eta; Q(k) \rangle = 1 + m k$. 
(See Section~\ref{sec:R21} for the definition of the map $Q: \mathbb{R} \rightarrow \mathfrak{g}$.)  
Therefore $(Q(k), \eta)$ belongs to $\mathcal{{\widetilde Z}}$ if and only if
\[
\eta = \eta(k, {\lambda}_1, \dots, {\lambda}_5) := \left( 1 + m k \right) \omega + 
{\lambda}_1 \eta^1 + {\lambda}_2 \eta^2 + {\lambda}_3 \eta^3 + {\lambda}_4 \eta^4 + {\lambda}_5 \eta^5,
\]
where ${\lambda}_1, \dots, {\lambda}_5 \in \mathbb{R}$. For simplicity, we identify ${\widetilde Z}$ with $G \times \mathbb{R}^6$ 
by means of the map
\[
(g: k, {\lambda}_1, \dots, {\lambda}_5) \in G \times \mathbb{R}^6 \rightarrow 
(g, Q(k), \eta(k, {\lambda}_1, \dots, {\lambda}_5)) \in {\widetilde Z}.
\]
Thus the Liouville form on ${\widetilde Z}$ is given by
\[
\psi = \left( 1 + m k \right) \omega + {\lambda}_1 \eta^1 + {\lambda}_2 \eta^2 + {\lambda}_3 \eta^3 + {\lambda}_4 \eta^4 + {\lambda}_5 \eta^5.
\]

From the structure equations~\eqref{structure}, we find that 
\begin{align*}
\Psi \equiv& \ 
m \pi \wedge \omega - \left( 1 + m \kappa \right) \eta^2 \wedge \omega
+ \sum d{\lambda}_{\alpha} \wedge \eta^\alpha
- {\lambda}_1 {\boldsymbol\pi} \wedge \omega 
+ {\lambda}_2 \left( \kappa \eta^3 - \eta^1 \right) \wedge \omega\\
&
+ {\lambda}_3 \left( 2 \eta^2 - \kappa \eta^4 \right) \wedge \omega 
+ {\lambda}_4 \left( \kappa \eta^5 - \eta^3 \right) \wedge \omega
+ {\lambda}_5 \eta^4 \wedge \omega,
\end{align*}
where $\Psi := d\psi$ and where $\equiv$ denotes equality modulo 
span$\left( \{ \eta^{\alpha} \wedge \eta^{\beta} \}_{\alpha, \beta = 1, \dots, 5} \right)$. 
Let
\[
\left( \partial_{\omega}, \partial_{\eta^1}, \dots \partial_{\eta^5}, 
\partial_{{\lambda}_1}, \dots, \partial_{{\lambda}_5}, \partial_{\pi} \right)
\]
denote the parallelization of ${\widetilde Z}$ dual to the coframing 
\[
\left( \omega, \eta^1, \dots \eta^5, d{\lambda}_1, \dots, d{\lambda}_5, \pi \right)
\]
We then have
\[
i_{\partial_{{\lambda}_i}} \Psi = \eta^i, \quad i = 1, \dots, 5,
\]
along with
\[
i_{\partial_{\omega}} \Psi \equiv - \alpha, \quad
i_{\partial_{\pi}} \Psi \equiv - \beta, \quad
i_{\partial_{\eta^i}} \Psi \equiv - \beta_i, \quad i = 1, \dots, 5,
\]
where
\begin{subequations}
\begin{align}
\alpha &= \left( m - {\lambda}_1 \right) dk,
\label{eq:alpha}\\
\beta &= \left( m - {\lambda}_1 \right) \omega,
\label{eq:beta}\\
\beta_1 &= d{\lambda}_1 + {\lambda}_2 \omega,\\
\beta_2 &= d{\lambda}_2 + \left( 1 + m k - 2 {\lambda}_3 \right) \omega,\\
\beta_3 &= d{\lambda}_3 + \left( {\lambda}_4 - k {\lambda}_2 \right) \omega,\\
\beta_4 &= d{\lambda}_4 + \left( k {\lambda}_3 - {\lambda}_5 - k \right) \omega,\\
\beta_5 &= d{\lambda}_5 - k {\lambda}_4 \omega,
\end{align}
\end{subequations}
and where $\equiv$ denotes equality modulo span$\left( \eta^1, \dots \eta^5 \right)$. 
From these equations, we deduce that the Cartan system $(\mathcal{C}(\Psi), \omega)$ is generated by 
the differential $1$-forms $(\eta^1, \dots, \eta^5, \alpha, \beta, \beta_1, \dots, \beta_5)$. 

\begin{thm}
\label{Y}
The momentum space, $Y$, is the $9$-dimensional sub-manifold of ${\widetilde Z}$ defined by the equations
\[
{\lambda}_2 = {\lambda}_1 - m = {\lambda}_3 - \frac{1}{2} \left( 1 + m k \right) = 0.
\]
and the Euler-Lagrange system $(\mathcal{E}, \omega)$ is the Pfaffian differential system on $Y$ 
with independence condition $\omega$ generated by the linear differential forms 
$(\eta^1, \dots, \eta^5, \sigma_1, \sigma_2, \sigma_3)$ where 
\[
\sigma_1 = \frac{m}{2} dk + {\lambda}_4 \omega,\quad
\sigma_2 = d{\lambda}_4 - \left( {\lambda}_5 + \frac{1}{2} k \left( 1 - m k \right) \right) \omega, \quad
\sigma_3 = d{\lambda}_5 - k {\lambda}_4 \omega.
\]
\end{thm}

\begin{proof}
We let $V_1 \subset T({\widetilde Z})$ be the sub-variety of $1$-dimensional integral elements of the Cartan system 
$(\mathcal{C}(\Psi), \omega)$ and denote by ${\widetilde Z}_1 \subset {\widetilde Z}$ the projection of $V_1$ under the bundle map 
$T({\widetilde Z}) \rightarrow {\widetilde Z}$. From equations~\eqref{eq:alpha} and~\eqref{eq:beta} we then deduce that ${\widetilde Z}_1$ 
is the submanifold defined by ${\lambda}_1 = m$. Denote by $\mathcal{C}(\Psi)_1$ the restriction to ${\widetilde Z}_1$ of the 
Cartan system. Then $\mathcal{C}(\Psi)_1$ is generated by the linear differential forms 
$(\eta^1, \dots, \eta^5, {\lambda}_2 \omega, \beta_2, \dots, \beta_5)$. 
We then consider the sub-variety $V_2 \subset T({\widetilde Z}_1)$ consisting of integral elements of 
$( \mathcal{C}(\Psi)_1, \omega )$ and let ${\widetilde Z}_2 \subset {\widetilde Z}_1$ denote the projection of $V_2$. 
We therefore have that ${\widetilde Z}_2$ is the sub-manifold of ${\widetilde Z}_1$ defined by ${\lambda}_2 = 0$. 
Denote by $\mathcal{C}(\Psi)_2$ the restriction to ${\widetilde Z}_2$ of $\mathcal{C}(\Psi)_1$.
Then $\mathcal{C}(\Psi)_2$ is generated by the linear differential forms 
$(\eta^1, \dots, \eta^5, \left( 1 + m k - 2 {\lambda}_3 \right) \omega, \beta_3, \beta_4, \beta_5)$. 
We proceed as above and let $V_3 \subset T({\widetilde Z}_2)$ be the sub-variety of integral elements 
of $( \mathcal{C}(\Psi)_2, \omega )$ and define ${\widetilde Z}_3 \subset {\widetilde Z}_2$ to be the image of $V_3$ 
under the projection $T({\widetilde Z}_2) \rightarrow {\widetilde Z}_2$. It follows that ${\widetilde Z}_3$ is the sub-manifold of 
${\widetilde Z}_2$ defined by the equation ${\lambda}_3 = \frac{1}{2} \left( 1 + m k \right)$ and that the 
restriction $\mathcal{C}(\Psi)_3$ of $\mathcal{C}(\Psi)_2$ to ${\widetilde Z}_3$ is the Pfaffian differential 
system generated by $(\eta^1, \dots, \eta^5, \sigma_1, \sigma_2, \sigma_3)$. 
If we let $V_4 \subset T({\widetilde Z}_3)$ be the set of integral elements of $( \mathcal{C}(\Psi)_3, \omega )$ 
then the bundle map $V_4 \rightarrow {\widetilde Z}_3$ is surjective. Hence $Y = {\widetilde Z}_3$ and 
$( \mathcal{C}(\Psi)_3, \omega )$ is the reduced space of $( \mathcal{C}(\Psi), \omega )$.
\end{proof}

\begin{corollary}
The momentum space $Y$ associated with the Poincar{\'e} variational problem for 
isotropic curves in {\Mink} is the rank-$3$, affine sub-bundle 
$Y = G \times \mathcal{F} \subset T^*(G) \cong G \times \mathfrak{g}^*$ where 
$\mathcal{F} \subset \mathfrak{g}^*$ is defined by
\[
\mathcal{F} = \frac{1}{2} \left( \omega^1 + \omega^2_1 \right) + m \omega^1_2 + 
\mathrm{span} \left( \omega^2_1 - \omega^1, \omega^2, \omega^3 \right).
\]
The variational problem is non-degenerate, and the characteristic vector field 
takes the form 
\begin{equation}
\xi = \partial_{\omega} - \frac{2{\lambda}_4}{m} \partial_k - {\lambda}_4 \partial_{{\lambda}_3} + 
\left( {\lambda}_5 + \frac{1}{2} k \left( 1 - m k \right) \right) \partial_{{\lambda}_4} + 
k {\lambda}_4 \partial_{{\lambda}_5}.
\label{xi}
\end{equation}
\end{corollary}

\begin{proof}
It follows from the preceding theorem that the restriction of the Liouville to the momentum space takes the form
\begin{align}
\psi_Y &= \left( 1 + m k \right) \omega + m \eta^1 + \frac{1}{2} \left( 1 + m k \right) \eta^3 + {\lambda}_4 \eta^4 + {\lambda}_5 \eta^5 
\nonumber \\
&= \frac{1}{2} \left( \omega^1 + \omega^2_1 \right) + m \omega^1_2 + 
\frac{mk}{2} \left( \omega^2_1 - \omega^1 \right) + {\lambda}_4 \omega^2 + {\lambda}_5 \omega^3.
\label{liouville}
\end{align}
The form of $Y$ and $\mathcal{F}$ follow directly from this equation. 
The dimension of $Y$ equals $9$, and a straightforward calculation shows that
\[
\omega \wedge \left( \Psi_Y \right)^4 = 
- 12 \, m^2 \, \omega \wedge dk \wedge d{\lambda}_4 \wedge d{\lambda}_5 
\wedge {\eta}^1 \wedge {\eta}^2 \wedge {\eta}^3 \wedge {\eta}^4 \wedge {\eta}^5,
\]
which is nowhere vanishing. Hence the variational problem is non-degenerate. 
The form of the characteristic vector field follows from a direct calculation.
\end{proof}

\begin{remark}
Since the variational problem is non-degenerate, 
the Euler-Lagrange system $\mathcal{E}$ coincides with the Cartan system of $\Psi$. 
The characteristic line-distribution $\Xi \subset T(Y)$ of $\Psi$ is transverse to the 
independence condition $\omega$, and is generated by the characteristic vector field $\xi$.
\end{remark}

\begin{remark}
Using the explicit form of the Liouville form, 
we may identify $Y = G \times \mathbb{R}^3$, 
where $(k, {\lambda}_4, {\lambda}_5)$ serve as coordinates on $\mathbb{R}^3$. 
The explicit form for the characteristic vector field and the $1$-forms $\eta^i$ and $\omega$ then imply that 
the map $\mathcal{H}: Y \rightarrow \mathfrak{g}$ is given by 
\[
\mathcal{H}[\eta(k, {\lambda}_4, {\lambda}_5)] = 
\begin{pmatrix}
0&0&0&0\\
1&0&k&0\\
0&1&0&k\\
0&0&1&0
\end{pmatrix} \in \mathfrak{g}.
\]
\end{remark}

\medskip

A smooth map $\Gamma: (a, b) \rightarrow Y$ is an integral curve of the Euler-Lagrange system 
if and only if it satisfies
\begin{align*}
\Gamma^* \left( \eta^i \right) &= 0, &i = 1, \dots, 5,
\\
\Gamma^* \left( \sigma_i \right) &= 0, &i = 1, 2, 3,
\end{align*}
with the independence condition
\[
\Gamma^* \left( \omega \right) \neq 0.
\]
Without loss of generality, we may choose a parametrization of our integral curve such that
\[
\Gamma^* \left( \omega \right) = dt.
\]
In this case, we may write 
$\Gamma: t \in (a, b) \rightarrow \left( g(t), k(t), {\lambda}_4(t), {\lambda}_5(t) \right) \in Y = G \times \mathbb{R}^3$. 
From these relations, and the explicit form of the differential forms $\eta^i$ and $\sigma_i$, we deduce the following result. 

\begin{proposition}
The smooth map $\Gamma: (a, b) \rightarrow Y$, parametrized such that $\Gamma^* \left( \omega \right) = dt$, 
is an integral curve of the Euler-Lagrange system if and only if the real-valued functions 
$k(t), {\lambda}_4(t), {\lambda}_5(t)$ satisfy the relations
\begin{equation}
\frac{dk}{dt} = - \frac{2{\lambda}_4}{m}, \qquad
\frac{d{\lambda}_4}{dt} = \left( {\lambda}_5 + \frac{1}{2} k \left( 1 - m k \right) \right), \qquad
\frac{d{\lambda}_5}{dt} = k {\lambda}_4,
\label{eulerr21}
\end{equation}
and $g(t) \in G$ is a solution of
\begin{equation}
g(t)^{-1} \frac{dg(t)}{dt} = \mathcal{H}(\eta(t)) = 
\begin{pmatrix}
0&0&0&0\\
1&0&k(t)&0\\
0&1&0&k(t)\\
0&0&1&0
\end{pmatrix}
\label{gdiffeqn}
\end{equation}
\end{proposition}

\begin{remark}
Although, in the present case, the Lie group $G$ is not semi-simple, it is naturally embedded in 
SL$(4, \mathbb{R})$. Using the Killing form on $\mathfrak{sl}(4, \mathbb{R})$ we deduce 
that the Euler equation~\eqref{eulerr21} may be written in Lax form: 
\[
L^{\prime} = \left[ L, \mathcal{H} \right],
\]
where 
\[
L(k, {\lambda}_4, {\lambda}_5) = 
\begin{pmatrix}
0&0&0&0\\
\frac{1}{2} \left( 1 + m k \right) &-{\lambda}_4 &-{\lambda}_5 &0\\
0&\frac{1}{2} \left( 1 - m k \right) &0 &-{\lambda}_5 \\
-m &0 &\frac{1}{2} \left( 1 - m k \right) &{\lambda}_4
\end{pmatrix}.
\]
\end{remark}
\smallskip

\begin{proposition}
If $\Gamma: (a, b) \rightarrow Y$ is an integral curve of the Euler-Lagrange system, 
with $\Gamma^* \left( \omega \right) = dt$, then the curvature $k(t)$ satisfies the third 
order ordinary differential equation
\begin{equation}
m \frac{d^3k}{dt^3} - 3 m k \frac{dk}{dt} + \frac{dk}{dt} = 0.
\label{3ode}
\end{equation}
Conversely, any non-constant solution $k: (a, b) \rightarrow \mathbb{R}$ of this equation determines 
a parametrized integral curve of the Euler-Lagrange system, unique up to the action of $\mathbb{E}(2, 1)$. 
\end{proposition}
\begin{proof}
Equation~\eqref{3ode} for $k(t)$ follows directly from~\eqref{eulerr21}. 
Conversely, given a solution of~\eqref{3ode}, 
we can uniquely reconstruct $\lambda_4(t), \lambda_5(t)$ from \eqref{eulerr21}, 
and $g(t)$ is determined, up to initial conditions, by \eqref{gdiffeqn}.
\end{proof}

\begin{remark}
Under the identification $\mathfrak{g}^* \cong \Mink \oplus \Mink$ introduced in Section~\ref{sec:coadj}, 
the Liouville form~\eqref{liouville} maps to $(p, v) \in \Mink \oplus \Mink$ where 
\[
p =  -  {\lambda}_5 e_1 +  {\lambda}_4 e_2 - \frac{1}{2} ( 1 - m k ) e_3, \quad
v =  - \frac{1}{2} ( 1 + m k ) e_1 +  m e_3
\]
and $\{ e_i \}$ is the standard basis of {\Mink}.
The Casimir operators therefore take the form 
\begin{equation}
C_1 = \Vert p \Vert^2 = {\lambda}_4^2 - {\lambda}_5 ( 1 - m k ), \quad
C_2 = \langle p, v \rangle = m {\lambda}_5 - \frac{1}{4} \left( 1 - m^2 k^2 \right),
\label{casimirs}
\end{equation}
and are constant along integral curves of the Euler-Lagrange system. 
\end{remark}

The explicit form of the Casimir operators implies the following result. 

\begin{proposition}
If $\Gamma: (a, b) \rightarrow Y$, parametrized such that $\Gamma^* \left( \omega \right) = dt$, 
is an integral curve of the Euler-Lagrange system then the curvature $k(t)$ satisfies the first
order ordinary differential equation
\[
\left( \frac{dk}{dt} \right)^2 = k^3 - \frac{1}{m} k^2 - \frac{1}{m^2} \left( 4 C_2 + 1 \right) k + 
\frac{1}{m^3} \left( 4 m C_1 + 4 C_2 + 1 \right).
\]
\end{proposition}

\begin{remark}
Letting $h(t) := \frac{1}{4} \left( k - \frac{1}{3m} \right)$, we deduce that $h$ satisfies
\[
\left( \frac{dh}{dt} \right)^2 = 4 h^3 - g_2 h - g_3,
\]
where 
\[
g_2 = \frac{1}{m^2} \left( C_2 + \frac{1}{3} \right), \qquad
g_3 = \frac{1}{m^3} \left( \frac{m C_1}{4} + \frac{C_2}{6} + \frac{1}{27} \right).
\]
Hence the curvature $k$ and the functions ${\lambda}_4, {\lambda}_5$ corresponding to any solution of 
the Euler-Lagrange system may be expressed in terms of Weierstrass 
elliptic functions with invariants $g_2, g_3$. 
\end{remark}

\begin{remark}
A short calculation using the explicit form of $\eta = \psi_Y$ given in~\eqref{liouville} 
and the coadjoint action of $\mathfrak{g}$ on $\mathfrak{g}^*$, which can be derived from~\eqref{coadjointaction},  
shows that in the present case the linearised phase-portrait $\Pi(\eta) = F(\eta) \cap O(\eta)$ is one-dimensional, 
and is spanned by the vector 
$ad^*[\mathcal{H}(\eta)] \eta = {\lambda}_4 \left( \omega^2_1 - \omega^1 \right) - 
\left( {\lambda}_5 + \frac{1}{2} k \left( 1 - m k \right) \right) \omega^2 - 
k {\lambda}_4 \omega^3$. Hence the regular parts of the phase portraits are one-dimensional in the current problem. 
We now introduce a more general class of variational problems for which this is the case.
\end{remark}

\section{Coisotropic variational problems}
\label{sec:coisotropic}
Consider a smooth manifold $M$ equipped with an exterior differential 2-form $\Psi$.
The kernel of $\Psi_x$ will be denoted by $N(\Psi)_x$.
Suppose that a Lie group $G$ acts on $M$.
Denote by $A^{\sharp}$ the fundamental vector field on $M$ corresponding to $A \in \mathfrak{g}$ and,
for each $x \in M$, let $\mathfrak{g}^{\sharp}(M)_x \subset T_x(M)$ be the vector subspace 
$\{ A^{\sharp}_x: A \in \mathfrak{g} \}$. 
We denote by $\mathfrak{g}^{\sharp}(M)^{\perp}_x$ the \emph{polar space\/} 
of $\mathfrak{g}^{\sharp}(M)_x$ with respect to $\Psi_x$: 
\[
\mathfrak{g}^{\sharp}(M)^{\perp}_x:= \left\{ 
v \in T_x(M): \Psi_x \left( v, A^{\sharp} \right) = 0, \forall A^{\sharp} \in \mathfrak{g}^{\sharp}(M)_x 
\right\}.
\]

\begin{definition}
The action of $G$ on $M$ is \emph{coisotropic\/} with respect to $\Psi$ if
\[
\mathfrak{g}^{\sharp}(M)^{\perp}_x \subset \mathfrak{g}^{\sharp}(M)_x+N(\Psi)_x, \quad \forall x \in M.
\]
\end{definition}

\begin{remark}
The notion of a coisotropic action arises naturally when studying collective complete integrability of
Hamiltonian systems (see \cite{GuilleminSternberg3,GuilleminSternberg1,Podesta}).
\end{remark}

\begin{definition}
An invariant Lagrangian $L: \mathbb{A} \to \mathbb{R}$ is said to be \emph{coisotropic\/}
if it is non-degenerate and if the action of $G$ on the regular part of the momentum space $Y$
is coisotropic with respect to $\Psi_{Y}$, the restriction of the canonical two-form $\Psi$ to $Y$.
\label{def:coisotropic}
\end{definition}

(Recall that we are now using the notation $Y$, $\mathcal{F}$ and $\mathcal{P}(\mu)$ to denote the
regular parts of the momentum space, phase space and phase portraits, respectively.)

\begin{proposition}
A non-degenerate invariant Lagrangian $L: \mathbb{A} \to \mathbb{R}$
is coisotropic if and only if the linearized phase portrait $\Pi(\eta)$
is spanned by $ad^*[\mathcal{H}(\eta)] \eta$, for every $\eta \in \mathcal{F}$.
\end{proposition}

\begin{proof}
Using the left-invariant trivialisation, 
we find that the polar space of $\mathfrak{g}^{\sharp}(Y)_{(g, \eta)}$ is given by
\[
\mathfrak{g}(\eta)^{\perp} := \mathfrak{g}^{\sharp}(Y)_{(g, \eta)}^{\perp} = 
\{ A + V \in \mathfrak{g} \oplus F(\eta): V = - ad^*(A) \eta \}.
\]

First, assume that $L$ is coisotropic, i.e.
\[
\mathfrak{g}(\eta)^{\perp} \subset \mathfrak{g} + \mathrm{span} \left( \mathcal{H}(\eta) - ad^*[\mathcal{H}(\eta)] \eta \right).
\]
Let $V \in \Pi(\eta)$. Then $V \in F(\eta)$ and there exists an $A \in \mathfrak{g}$ such that $V = - ad^*(A) \eta$.
Then $A + V$ belongs to $\mathfrak{g}(\eta)^{\perp}$ and hence $V$ must be a real multiple of $ad^*[\mathcal{H}(\eta)] \eta$.

Conversely, assume that $\Pi(\eta)$ is spanned by $ad^*[\mathcal{H}(\eta)] \eta$.
Given any element, $A + V$, of the polar space $\mathfrak{g}(\eta)^{\perp}$,
then $V \in \Pi(\eta)$ so there exists $s \in \mathbb{R}$ such that
$V = s \cdot ad^*[\mathcal{H}(\eta)] \eta$.
Therefore, we can write
\[
A + V = - s \left( \mathcal{H}(\eta) - ad^*[\mathcal{H}(\eta)] \eta \right) + \left( A +s \mathcal{H}(\eta) \right).
\]
Since $A + s \mathcal{H}(\eta)$ is an element of $\mathfrak{g}$, it follows that
$A + V \in \mathfrak{g} + \mathrm{span}\left( \mathcal{H}(\eta) - ad^*[\mathcal{H}(\eta)]\eta \right)$.
Therefore $\mathfrak{g}(\eta)^{\perp} \subset \mathfrak{g} +\mathrm{span}\left(\mathcal{H}(\eta)-ad^*[\mathcal{H}(\eta)]\eta\right)$, 
as required.
\end{proof}

\begin{remark}
Note that if $\eta \in \mathcal{F}$ then $\Pi(\eta)$ is $1$-dimensional
and the map $\rho(\eta): \mathrm{Ann}\!\left( F(\eta) \right) \to O(\eta)$ defined in~\eqref{rhodef} is injective.
\end{remark}

\begin{proposition}
\label{prop:rank}
Let $L: \mathbb{A} \to \mathbb{R}$ be an invariant coisotropic
Lagrangian and let $Y = G \times \mathcal{F}$ be the corresponding momentum space. 
Suppose that $\mathcal{F}$ is non empty, we then have:
\begin{itemize}
\item
$\mathrm{dim}(Y) = \mathrm{dim}(G) + \mathrm{rank}(G) + 1$;
\item the regular part of the phase space, $\mathcal{F}$, 
intersects the coadjoint orbits transversally;
\item The regular parts $\mathcal{P}(\mu)$ of the phase portraits are smooth and $1$-dimensional;
\item Every $\eta \in \mathcal{F}$ is a regular element of $\mathfrak{g}^*$. 
In particular, the isotropy group $G_{\eta}$ and the isotropy algebra
$\mathfrak{g}_{\eta}$ are Abelian.
\end{itemize}
\end{proposition}

\begin{proof}
For each $\eta \in \mathcal{F}$, let $k(\eta)$ be the
dimension of the isotropy Lie algebra $\mathfrak{g}_{\eta}$. 
Note that
\[
\mathrm{dim}(F(\eta)\cap O(\eta)) = 1, \quad
\mathrm{dim}(O(\eta)) = \mathrm{dim}(G) - k(\eta), \quad
\mathrm{dim}\left(\mathrm{Ann}\!\left(F(\eta)\right)\right) = \mathrm{dim}(G) - \mathrm{dim}(F(\eta)).
\]
We then have
\[
\mathrm{dim}(F(\eta)) + \mathrm{dim}(O(\eta)) - 1 \leq \mathrm{dim}(G),
\]
which in turn implies that $\mathrm{dim}(\mathcal{F}) \leq k(\eta) + 1$. 
On the other hand, from the injectivity of the map 
$\rho(\eta): \mathrm{Ann}\!\left(F(\eta)\right) \to O(\eta)$ it follows that
$k(\eta) \leq \mathrm{dim}(\mathcal{F})$. 
Therefore we have
\[
k(\eta) \leq \mathrm{dim}(\mathcal{F}) \leq k(\eta) + 1.
\]
Notice that $\mathrm{dim}(G) + k(\eta)$ is even and that
$\mathrm{dim}(Y) = \mathrm{dim}(G) + \mathrm{dim}(\mathcal{F})$ is odd. 
Thus, we must have $k(\eta) + 1 = \mathrm{dim}(\mathcal{F})$. 
In particular, $k(\eta) = k$ is constant and
\[
\mathrm{dim}(Y) = \mathrm{dim}(G) + k + 1.
\]
This implies that
\[
\mathrm{dim}(\mathfrak{g}) = \mathrm{dim}(F(\eta)) + \mathrm{dim}(O(\eta)) - 1.
\]
Thus $\mathcal{F}$ intersects the coadjoint orbits transversally. 
Since $\mathrm{dim}(F(\eta)) \cap \mathcal{O}(\eta)) = 1$, 
it follows that $\mathcal{P}(\eta) = \mathcal{F} \cap \mathcal{O}(\eta)$ 
is a smooth curve such that $T_{\eta}[\mathcal{P}(\eta)] = \Pi(\eta)$. 
Moreover, from the transversality condition, it follows that 
$\mathcal{F}$ cannot be contained in the set $\mathfrak{g}^*_s$ 
of the singular element of $\mathfrak{g}^*$.
Thus, $\mathcal{F} \cap \mathfrak{g}^*_r$ is non-empty. 
Therefore, there exists an $\eta \in \mathcal{F}$ such that 
$k = k(\eta) = \mathrm{rank}(G)$. 
This gives the required result.
\end{proof}

\begin{remark}
The regular part of the phase space, $\mathcal{F}$, 
is foliated by the nowhere vanishing vector field $\Phi_{\xi}$ 
and the leaves are the phase portraits. 
Furthermore, if $X \subset \mathcal{F}$ is a local section of such a foliation 
then $X$ is also a local section of the coadjoint representation.
\end{remark}

\begin{definition}
Let $L: \mathbb{A} \to \mathbb{R}$ be a coisotropic Lagrangian.
The \emph{moment map\/} $J: Y \to \mathfrak{g}^*$ of the Hamiltonian action of $G$ on $Y$ is defined by
\[
J(g, \eta) = Ad^*(g) \eta, \quad \forall (g, \eta) \in Y.
\]
\end{definition}

\begin{proposition}
Let $L: \mathbb{A} \to \mathbb{R}$ be a coisotropic Lagrangian. 
Then
\begin{itemize}
\item $J(Y) \subset \mathfrak{g}_r^*$;
\item $J: Y \rightarrow \mathfrak{g}^*$ is a submersion; 
\item $J^{-1}(\mu)$ is a $(k+1)$-dimensional submanifold of $Y$ such that
\[
T_{(g, \eta)}[J^{-1}(\mu)] = \ker[dJ|_{(g, \eta)}] = \mathrm{span}[\xi|_{(g, \eta)}] + \mathfrak{g}_{\eta}.
\]
so that the characteristic vector field $\left. \xi \right|_{J^{-1}(\mu)}$ is tangent of $J^{-1}(\mu)$, 
and $G_{\mu}$ acts freely and properly on $J^{-1}(\mu)$. 
\item $Y_{\mu} := J^{-1}(\mu)/G_{\mu}$ is a one-dimensional manifold and $J^{-1}(\mu) \rightarrow Y_{\mu}$ is a 
principal $G_{\mu}$ bundle.
\item $Y_{\mu} \cong \mathcal{P}(\mu)$, the phase portrait.
\end{itemize}
\end{proposition}

\begin{proof}
From Proposition~\ref{prop:rank} we know that each $\eta \in \mathcal{F}$ is an element of $\mathfrak{g}^*_r$, 
and hence $J(g, \eta) \in \mathfrak{g}^*_r$ for all $(g, \eta) \in Y$. 
The differential of the moment map is given by the formula
\begin{equation}
\langle dJ|_{(g, \eta)}(A+V); B \rangle =
\langle ad^*(A) \eta + V; B \rangle, \quad
\forall A+V \in \mathfrak{g} \oplus F(\eta), 
\quad \forall B \in \mathfrak{g}.
\label{dJ}
\end{equation}
This implies that 
\[
\mathrm{Im}\!\left[ dJ|_{(g, \eta)} \right]= F(\eta) + O(\eta),
\]
for all $\eta \in \mathcal{F}$. 
Since $\mathcal{F}$ and $\mathcal{O}(\eta)$ intersect transversally, 
this implies that $J$ is a submersion. 
Therefore $J^{-1}(\mu)$ is a sub-manifold of $Y$ and the tangent space 
$T_{(g, \eta)}[J^{-1}(\mu)]$ is naturally isomorphic to $\ker[dJ|_{(g, \eta)}]$. 
From the formula~\eqref{dJ} and the fact that the Lagrangian is coisotropic, 
we deduce that
\[
\ker\!\left[ dJ|_{(g, \eta)} \right] = \mathrm{span}[\xi|_{(g, \eta)}] + \mathfrak{g}_{\eta}, 
\]
for all $\eta \in \mathcal{F}$, as required. 
Note that this relation implies that the characteristic vector field $\xi$ belongs to $\ker \left[ dJ \right]$, 
and therefore that $\left. \xi \right|_{(g, \eta)}$ is tangent to $J^{-1}(\mu)$. 
We shall denote the restriction of $\xi$ to the fibre $J^{-1}(\mu)$ by $\xi^{\mu}$. 

The isotropy group $G_{\mu}$ acts on $J^{-1}(\mu)$ by $(g, \eta) \mapsto (h g, \eta)$ 
for each $h \in G_{\mu}$. This action is clearly free and proper, 
so the quotient space $Y_{\mu} := J^{-1}(\mu)/G_{\mu}$ 
exists as a one-dimensional manifold. The map 
\[
\pi_{\mu}: (g, \eta) \in J^{-1}(\mu) \to [(g, \eta)] \in Y_{\mu}
\]
gives $J^{-1}(\mu)$ the structure of a 
principal fibre-bundle with structure group $G_{\mu}$. 
Moreover, the vector field $\xi^{\mu}$ is horizontal
with respect to the fibration $J^{-1}(\mu) \to Y_{\mu}$.

We also consider the fibration of $J^{-1}(\mu)$ over $\mathcal{P}(\mu)$ defined by 
\[
\widetilde{\pi}_{\mu}: (g, \eta) \in J^{-1}(\mu) \to \eta \in \mathcal{P}(\mu).
\]
The structure group is again the isotropy subgroup $G_{\mu}$. 
Furthermore, $\widetilde{\pi}_{\mu}$ is constant along the fibers of the
fibration $\pi_{\mu}$, and therefore descends
to a diffeomorphism of $Y_{\mu}$ onto $\mathcal{P}(\mu)$. 
\end{proof}

\begin{definition}
We adopt standard terminology, referring to $Y_{\mu}$ as the
\emph{Marsden-Weinstein reduction\/} of $Y$ at $\mu$, 
and to $\pi_{\mu}: J^{-1}(\mu) \to Y_{\mu}$ as the
\emph{Marsden-Weinstein fibration\/} at $\mu$. 
We may consider $\pi_{\mu}: J^{-1}(\mu) \to Y_{\mu}$ 
as a principal $G_{\mu}$ bundle over $Y_{\mu}$, 
where the right-action of $G_{\mu}$ on $J^{-1}(\mu)$ 
is given by $R_{h}(g, \eta)=(h g, \eta)$, 
for all $h \in G_{\mu}$, for $(g, \eta) \in J^{-1}(\mu)$.
\end{definition}

\begin{definition}
Let $\mu \in J(Y)$. The restriction of the Marsden-Weinstein
fibration $J^{-1}(\mu)$ to the reduced phase portrait
$\widetilde{\mathcal{P}}(\mu)$ is said to be the \emph{reduced Marsden-Weinstein fibration}. 
We shall denote this fibration by $\widetilde{\pi}_{\mu}: P^{\mu} \to \widetilde{\mathcal{P}}(\mu)$.
\end{definition}

\begin{remark}
The vector field $\Phi_{\xi}: \eta \mapsto -ad^*[\mathcal{H}(\eta)] \eta$ 
is tangent to the phase portraits. 
We denote by $\Phi_{\xi}^{\mu}$ the restriction of $\Phi_{\xi}$ to $\widetilde{\mathcal{P}}(\mu)$. 
Note that the vector fields $\xi^{\mu}$ and $\Phi_{\xi}^{\mu}$ 
are related by the fibration $\tilde{\pi}_{\mu}: P^{\mu} \to \widetilde{\mathcal{P}}(\mu)$.
\end{remark}

\begin{remark}
On the reduced phase portrait $\widetilde{\mathcal{P}}(\mu)$ 
there exists a unique nowhere vanishing $1$-form $\sigma^{\mu}$ 
such that $\sigma^{\mu}(\Phi_{\xi}^{\mu})=1$.
Take $\eta \in \widetilde{\mathcal{P}}(\mu)$, then the integral curve
$\phi_{\eta}: (\epsilon^-(\eta), \epsilon^+(\eta)) \to \mathfrak{g}^*$ 
is a maximal parametrization of $\widetilde{\mathcal{P}}(\mu)$ such that
$\phi_{\eta}^*(\sigma^{\mu})=dt$.
\end{remark}

\begin{definition}
On $P^{\mu}$ we consider the $\mathfrak{g}_{\mu}$-valued $1$-form 
$\theta^{\mu}$ defined by 
\[
\left. \theta^{\mu} \right|_{(g, \eta)} := \mathrm{Ad}(g) \left( \Theta - \mathcal{H} \sigma^{\mu} \right).
\] 
This defines a connection on the reduced Marsden-Weinstein fibration 
$P^{\mu} \to \widetilde{\mathcal{P}}(\mu)$. 
We call $\theta^{\mu}$ the \emph{canonical connection\/} 
of the reduced Marsden-Weinstein fibration $P^{\mu} \to \widetilde{\mathcal{P}}(\mu)$.
\end{definition}

\subsection{Isotropic curves in {\Mink}}
\label{sec:parametrisation}
In the case of our problem for isotropic curves in {\Mink}, we have defined a map 
\begin{equation}
\mathbb{R}^3 \hookrightarrow \mathfrak{g}^* \cong \Mink \oplus \Mink, \quad 
y = \left( k, {\lambda}_4, {\lambda}_5 \right) \mapsto \left( p(y), v(y) \right),
\label{embed}
\end{equation}
where 
\[
p = \begin{pmatrix} -  {\lambda}_5\\ {\lambda}_4\\ - \frac{1}{2} ( 1 - m k ) \end{pmatrix}, \qquad
v = \begin{pmatrix} - \frac{1}{2} ( 1 + m k )\\ 0\\ m \end{pmatrix}.
\]
Given the form~\eqref{xi} of the characteristic vector field $\xi$, 
we see that the regular part of the phase space $\mathcal{F}$ 
is given by the complement of the set of points with
\begin{equation}
{\lambda}_4 = 0, \quad
{\lambda}_5 + \frac{k}{2} \left( 1 - m k \right) = 0.
\label{singular}
\end{equation}

To show that the action of $\mathbb{E}(2, 1)$ on the regular part of the momentum space, $Y$, 
is coisotropic, we consider a general vector field on $Y$: 
\[
Z^1 \frac{\partial}{\partial k} + 
Z^2 \frac{\partial}{\partial {\lambda}_4} + 
Z^4 \frac{\partial}{\partial {\lambda}_5} + 
X^i \frac{\partial}{\partial \omega^i} + 
Y^1 \frac{\partial}{\partial \omega^1_1} + 
Y^2 \frac{\partial}{\partial \omega^2_1} + 
Y^3 \frac{\partial}{\partial \omega^1_2}.
\]
This vector field lies in $\left. \mathfrak{g}(Y)^{\perp} \right|_{(g, k, {\lambda}_4, {\lambda}_5)}$ if and only if 
\begin{subequations}
\begin{align}
\frac{m}{2} Z^1 &= - \frac{1}{2} \left( 1 - m k \right) Y^1 - {\lambda}_4 Y^2,\\
Z^2 &= \frac{1}{2} \left( 1 - m k \right) Y^3 + {\lambda}_5 Y^2,\\
Z^3 &= {\lambda}_4 Y^3 - {\lambda}_5 Y^1
\end{align}
\label{star}
\end{subequations}
and
\begin{subequations}
\begin{align}
- \frac{1}{2} \left( 1 - m k \right) X^1 + {\lambda}_5 X^3 - m Y^3 + \frac{1}{2} \left( 1 + m k \right) Y^2 &= 0,\\
- {\lambda}_4 X^1 - {\lambda}_5 X^2 - \frac{m}{2} Z^1 - \frac{1}{2} \left( 1 + m k \right) Y^1 &= 0,\\
- \frac{1}{2} \left( 1 - m k \right) X^2 - {\lambda}_1 X^3 + m Y^1 &= 0.
\end{align}
\label{bigstar}
\end{subequations}
If $(k, {\lambda}_4, {\lambda}_5)$ lies in the regular part of $\mathcal{F}$ then conditions~\eqref{star} and~\eqref{bigstar} are 
linearly independent. Therefore in this case
\[
\dim \mathfrak{g}(Y)^{\perp} = 3.
\]
It is easily checked that the characteristic vector field $\xi$ belongs to $\mathfrak{g}(Y)^{\perp}$, 
as do the following vector fields:
\begin{align*}
S_1 &= - {\lambda}_4 \frac{\partial}{\partial \omega^1_1} + 
\frac{1}{2} \left( 1 - m k \right) \frac{\partial}{\partial \omega^2_1} - 
{\lambda}_5 \frac{\partial}{\partial \omega^1_2} - m \frac{\partial}{\partial \omega^3},\\
S_2 &= {\lambda}_5 \frac{\partial}{\partial \omega^1} - 
{\lambda}_4 \frac{\partial}{\partial \omega^2} + 
\frac{1}{2} \left( 1 - m k \right) \frac{\partial}{\partial \omega^3}.
\end{align*}
Hence we have that 
\[
\mathfrak{g}(Y)^{\perp} = \mathrm{span} ( \xi, S_1, S_2 ) \subset \mathrm{span} ( \xi ) \oplus \mathfrak{g}.
\]
Hence the action is coisotropic.

\smallskip\noindent{\textbf{The momentum map and the basic invariants}}: 
From the form of the coadjoint action of $\mathbb{E}(2, 1)$ given earlier, 
we deduce that the moment map takes the form 
\[
J(g; y) = \left( A p(y), A v(y) - \left( A p(y) \right) \times Q \right),
\]
for all $g = \left( Q, A \right) \in \mathbb{E}(2, 1)$. Note that $J(Y) \subseteq \mathfrak{g}^*_r$. 
The basic invariants, which correspond to constants of motion of the system, are the Casimir operators 
\begin{align*}
C_1(g, y) := \Vert p(y) \Vert^2 = {\lambda}_4^2 - {\lambda}_5 \left( 1 - m k \right),\\
C_2(g, y) := \langle p(y), v(y) \rangle = m {\lambda}_5 - \frac{1}{4} \left( 1 - m^2 k^2 \right).
\end{align*}
If we choose $\mu = \left( \mathfrak{m}_1, \mathfrak{m}_2 \right) \in \mathfrak{g}_r^*$, 
where $\mathfrak{m}_1, \mathfrak{m}_1 \in {\Mink}$, then $J^{-1}(\mu)$ is the set
\begin{subequations}
\begin{align}
{\lambda}_4^2 &= \frac{m^2}{4} k^3 - \frac{m}{4} k^2 - \left( \frac{1}{4} + \langle \mathfrak{m}_1, \mathfrak{m}_2 \rangle \right) k + 
\left( \Vert \mathfrak{m}_1 \Vert^2 + \frac{1}{4m} + \frac{\langle \mathfrak{m}_1, \mathfrak{m}_2 \rangle}{4} \right),\\
{\lambda}_5 &= \frac{1}{m} \left( \frac{1}{4} \left( 1 - m^2 k^2 \right) + \langle \mathfrak{m}_1, \mathfrak{m}_2 \rangle \right).
\end{align}
\label{p4p5}
\end{subequations}
If we perform the substitution 
\begin{equation}
k = \left( \frac{4}{m} \right)^{2/3} \chi + \frac{1}{3m},
\label{ksub}
\end{equation}
then the first equation becomes the cubic relation
\begin{equation}
{\lambda}_4^2 = 4 \chi^3 - g_2 \chi - g_3,
\label{poly}
\end{equation}
where we have defined the \lq\lq modified Casimirs\rq\rq\ 
\begin{align*}
g_2(\mathfrak{m}_1, \mathfrak{m}_2) &= \left( \frac{4}{m} \right)^{2/3} 
\left( \frac{1}{3} + \langle \mathfrak{m}_1, \mathfrak{m}_2 \rangle \right),\\
g_3(\mathfrak{m}_1, \mathfrak{m}_2) &= - \left( 
\Vert \mathfrak{m}_1 \Vert^2 + \frac{2}{3m}\langle \mathfrak{m}_1, \mathfrak{m}_2 \rangle + \frac{4}{27m} \right).
\end{align*}

\noindent{\textbf{Parametrization of the phase portraits}}: 
We define the discriminant of the cubic polynomial appearing in equation~\eqref{poly}
\[
D(m_1, m_2) = 27 g_3^2 - g_2^3.
\]
There are two non-degenerate cases that we must consider.

\smallskip
\noindent{\textbf{Case I}: $D(m_1, m_2) > 0$}: 
In this case the cubic polynomial has one real root and two complex-conjugate roots. 
We may parametrise the curve 
by taking $\chi(t) = \wp(t; g_2, g_3)$ with $t \in \left( 0, 2 \omega_1 \right)$ where 
$\omega_1, \omega_2$ and $\omega_3 = \frac{1}{2} \left( \omega_1 + \omega_2 \right)$ 
are the half-periods of the $\wp$ function. From~\eqref{ksub}, we then find 
that $k$ may be written in terms of elliptic functions and, solving~\eqref{p4p5} 
then gives ${\lambda}_4$ and ${\lambda}_5$ in terms of elliptic functions.

\smallskip
\noindent{\textbf{Case II}: $D(m_1, m_2) < 0$}: 
In this case the cubic polynomial has 
three distinct real roots, and the curve \eqref{poly} has two disjoint components. 
The \lq\lq compact\rq\rq\ component may be parametrized by 
$\chi(t) = \wp_3(t; g_2, g_3) := \wp(t + \omega_3; g_2, g_3)$ 
with $t \in \mathbb{R}$. This solution is periodic, with period $2\omega_1$. 
The \lq\lq unbounded\rq\rq\ component of the curve is parametrized by 
$\chi(t) = \wp(t; g_2, g_3), t \in \left( 0, 2 \omega_1 \right)$. 

\smallskip

In the degenerate cases where the discriminant vanishes $D=0$, 
the cubic is singular, and our curve is rational. 
In this case, the $\wp$ and $\wp_3$ functions degenerate into elementary functions. 

\section{Integrability by quadratures}

\begin{proposition}
The integral curves of the characteristic vector
field $\xi$ with momentum $\mu \in J(Y)$ are the horizontal
curves of the canonical connection $\theta^{\mu}$ on $P^{\mu}$.
\end{proposition}

\begin{proof}
Consider a horizontal curve $\Gamma: (a, b) \to P^{\mu}$ of the canonical connection. 
We write $\Gamma(t)=(h(t), \eta(t))$, where 
$h: (a, b) \to G$ and $\eta: (a, b) \to \mathcal{P}(\mu)$ are smooth curves. 
Without loss of generality we can assume that $\eta^*(\sigma^{\mu}) = dt$, so that
\[
\eta'(t) = \left. \Phi^{\mu}_{\xi} \right|_{\eta(t)} = - ad^*[\mathcal{H}(\eta(t))] \eta(t), 
\quad \forall t \in (a, b).
\]
Since $\Gamma$ is horizontal, we then have
\[
0 = \Gamma^*(\theta^{\mu}) = g(t) (h^{-1}(t) h^{\prime}(t) - \mathcal{H}[\eta(t)]) g(t)^{-1}\, dt.
\]
This implies that $\Gamma'(t) = \xi|_{\Gamma(t)}$, 
for all $t \in (a, b)$. 

Conversely, if $\Gamma: (a, b) \to Y$ is an integral curve 
of $\xi$ with momentum $\mu$ then $\Gamma(t) \in P^{\mu}$, 
for all $t \in (a, b)$. 
Furthermore, we know that
\[
\Gamma'(t) = \mathcal{H}(\eta(t)) - ad^*[\mathcal{H}(\eta(t))] \eta(t).
\]
Thus $h^{-1}(t) h'(t) = \mathcal{H}[\eta(t)]$ and hence $\Gamma^*[\theta^{\mu}] = 0$.
\end{proof}
\vskip .2cm

\noindent{}Since the structure group of the Marsden-Weinstein
fibrations is Abelian and the base manifolds are
$1$-dimensional, the horizontal curves of the canonical
connection can be found by a single quadrature. The explicit
integration of the horizontal curves requires four steps:
\begin{itemize}
\item Step one: take a smooth parametrization 
$\eta: (a, b) \to \widetilde{\mathcal{P}}(\mu)$ of the phase portrait.
\item Step two: compute $v^{\mu}: (a, b) \to \mathbb{R}$ such that
$\eta^* (\sigma^{\mu}) = v^{\mu} dt$.
\item Step three: take any map $g: (a, b) \to G$ such that $(g^{-1}, \eta): (a, b) \to Y$
is a cross-section of the reduced Marsden-Weinstein fibration
${\widetilde\pi}_{\mu}: P^{\mu} \to \widetilde{\mathcal{P}}(\mu)$. 
This involves solving the equation
\[
Ad^*(g(t)^{-1}) \eta(t) = \mu
\] 
for $g(t)$.
\item Step four: compute the gauge transformation
\begin{equation}
h(t) = 
\mathrm{Exp}
\left[
\int_{t_0}^{t} \left( g^{-1}(u) \mathcal{H}[\eta(u)] g(u) v^{\mu}(u) + g^{-1}(u) g'(u) \right) du \right], \quad
\forall t \in (a, b).
\label{gaugetransformation}
\end{equation}
Note that the fact that $(g(t)^{-1}, \eta(t))$ is a section of 
${\widetilde\pi}_{\mu}: P^{\mu} \to \widetilde{\mathcal{P}}(\mu)$, 
along with the definitions of $\mathcal{H}[\eta]$ and $v^{\mu}$ imply that 
$g^{-1}(t) \mathcal{H}[\eta(t)] g(t) v^{\mu}(t) + g^{-1}(t) g'(t) \in \mathfrak{g}_{\mu}$ 
for all $t \in (a, b)$. Hence $h(t) \in G_{\mu}$, for all $t \in (a, b)$.
\end{itemize} 
\noindent{}\textbf{Conclusion}: 
the image of the curve $\Gamma: (a, b) \to Y$ defined by
\begin{equation}
\Gamma(t)=(h(t)g(t)^{-1}, \eta(t)), \quad \forall t \in (a, b)
\label{horizontalcurve}
\end{equation}
is contained in $P^{\mu}$ and 
$\Gamma$ is horizontal for the canonical connection $\theta^{\mu}$. 
Any horizontal curve of $P^{\mu}$ arises in this way.

\begin{remark}
In the case where the symmetry group $G$ is a classical matrix group, 
with $G \subset \mathrm{Aut}(V)$ for $V$ a finite-dimensional vector space, 
and the Lagrangian is a polynomial function, then we have:
\begin{itemize}
\item The phase space $\mathcal{F}$ is an algebraic subset of $\mathfrak{g}^*$.
\item The generic coadjoint orbit is defined by polynomial
equations $F_j(\eta) = 0, j = 1, \dots, k$, where $k$ is the rank of
$\mathfrak{g}$ and where $F_1, \dots, F_k$ is a basis of the $Ad^*$-invariant
polynomial functions $\mathfrak{g}^* \to \mathbb{R}$.
\item The phase portraits are real-algebraic curves.
\end{itemize}
In the simplest cases the phase portraits are rational or elliptic
curves, so they can be easily parametrized by means of
elementary or elliptic functions (see the examples considered in
refs.~\cite{BryantGriffiths,Griffiths1,Hsu,Musso1,Musso2,Musso3,MussoNicolodi}).
The third step in the construction above can be treated in a rather easy way if $\mathfrak{g}$ is
semi-simple and if the momentum $\mu$ is a regular semisimple
element of $\mathfrak{g}^*$. In this case the construction of a
cross-section of the Marsden-Weinstein fibration is a
linear-algebra problem involving the structure of the Cartan 
subalgebras of $\mathfrak{g}$.
\end{remark}

\begin{definition}
Consider $\mu \in J(Y)$. 
We say that $\mu$ is a \emph{complete momentum\/} if 
$\epsilon^{-}(\eta)=- \infty$ and $\epsilon^{+}(\eta)= \infty$ 
for some (and hence for all) $\eta \in \widetilde{\mathcal{P}}(\mu)$.
\end{definition}

\begin{proposition}
If $\mu \in J(Y)$ is a complete momentum, then the connected components of the reduced Marsden-Weinstein
fibration $P^{\mu} \to \widetilde{\mathcal{P}}(\mu)$ are Euclidean cylinders and $\xi^{\mu}$ is a linear vector field.
\end{proposition}

\begin{proof}
Let $Q(\mu)$ be a connected component of $P^{\mu}$ 
and let $\eta: \mathbb{R} \to \widetilde{P}(\mu)$ 
be an integral curve of the phase flow 
that parametrizes the reduced phase portrait. 
Since $\mathbb{R}$ is contractible, 
$Q(\mu) \to \widetilde{P}(\mu)$ is a trivial fiber bundle. 
This implies that there exists a smooth map $g: \mathbb{R} \to G$ 
such that $(g^{-1}, \eta): \mathbb{R} \to G \times \widetilde{P}(\mu)$ 
is a cross-section of $Q(\mu) \to \widetilde{P}(\mu)$. 
Fix $(g_0, \eta_0) \in Q(\mu)$ and let $t_0 \in \mathbb{R}$ 
such that $\eta(t_0) = \eta_0$ and consider the curve 
$\Gamma_{(g_0, \eta_0)}: \mathbb{R} \to Q(\mu)$ defined by
\[
\Gamma_{(g_0, \eta_0)}(t) = (g_0 g(t_0)k(t), \eta(t)), 
\quad \forall t \in \mathbb{R},
\]
where
\[
k(t) = \mathrm{Exp}\left[ \int_{t_0}^{t}
\left( g(u)^{-1} \mathcal{H}[\eta(u)] g(u) + g(u)^{-1} g'(u) \right) du \right] g(t)^{-1}, 
\quad \forall t \in \mathbb{R}.
\]
Then, $\Gamma_{(g_0, \eta_0)}$ is the integral curve of $\xi^{\mu}$ with initial condition
$\Gamma_{(g_0, \eta_0)}(t_0)=(g_0, \eta_0)$. 
This shows that the restriction of the vector field $\xi^{\mu}$ to $Q(\mu)$ is complete. 
Now fix a basis $(e_1, \dots, e_k)$ of $\mathfrak{g}_{\mu}$ and
let $e^{\sharp}_1, \dots, e^{\sharp}_k$ denote the corresponding
fundamental vector fields on $Q(\mu)$. 
Then $\{ \xi^{\mu}, e^{\sharp}_1, \dots, e^{\sharp}_k \}$ is a set of
complete, linearly independent and commuting vector fields on $Q^{\mu}$. 
It is then a standard fact that $Q(\mu)$ is a $(k+1)$-dimensional cylinder, 
that is $Q(\mu) = \mathbb{R}^{k+1}/K$, where $K \subset \mathbb{R}^{k+1}$ 
is a subgroup of $\mathbb{R}^{k+1}$ generated over $\mathbb{Z}$ by $m \le k+1$ 
linearly independent vectors $a_1, \dots, a_m$:
\[
K = \left\{ \sum_{j=1}^{m}n_j a_j: n_j \in \mathbb{Z} \right\}.
\]
The vector fields $\{\xi^{\mu}, e^{\sharp}_1, \dots, e^{\sharp}_k\}$ 
are then the push-forward of linear vector fields 
$b_0, b_1, \dots, b_k$ on $\mathbb{R}^{k+1}$. 
This yields the required result.
\end{proof}

\begin{remark}
If $\widetilde{\mathcal{P}}(\mu)$ is compact and the isotropy subgroup $G_{\mu}$ is compact, 
then the connected components of the reduced Marsden-Weinstein fibration 
$P^{\mu}$ are $(k+1)$-dimensional tori.
\end{remark}

\subsection{Cross-sections of the Marsden-Weinstein fibration for isotropic curves in \Mink}

Finally, we show how the above integration procedure may be carried out in our example. 

Given $y = \left( k, {\lambda}_4, {\lambda}_5 \right) \in \mathbb{R}^3$, 
we have defined the vectors 
$\left( p, v \right) \in \Mink \oplus \Mink \cong \mathfrak{g}^*$. 
With respect to the standard basis $(e_1, e_2, e_3)$ for {\Mink} 
these take the form
\[
p = \sum_{i=1}^3 p^i e_i := - {\lambda}_5 e_1 + {\lambda}_4 e_2 - \frac{1}{2} \left( 1 - m k \right) e_3, \qquad
v = \sum_{i=1}^3 v^i e_i :=- \frac{1}{2} \left( 1 + m k \right) e_1 + m e_3.
\]
Letting $\mu = (\mathfrak{m}_1, \mathfrak{m}_2) \in \mathrm{Im}\, J \subseteq \mathfrak{g}^*_r$, 
we wish to construct the map $g: (a, b) \subset \mathbb{R} \rightarrow G$ 
with the property that $(g^{-1}, \eta)$ is a section of the reduced Marsden-Weinstein fibration 
${\widetilde\pi}_{\mu}: P^{\mu} \to \widetilde{\mathcal{P}}(\mu)$. 
We must consider separately the cases where the coadjoint orbit is of positive, 
negative or null type. 

\medskip\noindent{\textbf{Positive type}}: 
For an orbit of positive type, where $C_1 = \Vert p \Vert^2 > 0$, 
we may assume, 
up to the action of $G$ on $\mathcal{O}(\mu)$, 
that $\mu = (\mathfrak{m}_1, \mathfrak{m}_2)$ is in the standard form:
\[
\mathfrak{m}_1 = \sqrt{C_1} \begin{pmatrix} 0\\ 1\\ 0 \end{pmatrix}, \quad
\mathfrak{m}_2 = \frac{C_2}{\sqrt{C_1}} \begin{pmatrix} 0\\ 1\\ 0 \end{pmatrix},
\]
where $C_2 := \langle p, v \rangle$. 


We now wish to construct $g = (Q,  A) \in \mathbb{E}(2, 1)$ with the property that
\[
\eta = \left( p, v \right) = \mathrm{Ad}^*(g) \mu.
\]
Since $\Vert p \Vert^2 > 0$, we may define a the $\Mink$-valued map
\[
A_2 := \frac{p}{\sqrt{C_1}}= \frac{p}{\Vert p \Vert},
\]
with the property that $\Vert A_2 \Vert^2 = 1$. 
We can now complete $A_2$ to a frame field $\left( A_1, A_2, A_3 \right)$ 
by adding any $\Mink$-valued functions $A_1, A_3$ with the property that 
\[
\langle A_i, A_j \rangle = g_{ij}, \qquad i, j = 1, 2, 3,
\]
and we fix the orientation of this basis by the requirements that
\[
A_2 \times A_1 = A_1, \qquad
A_2 \times A_3 = - A_3, \qquad
A_3 \times A_1 = A_2.
\]

More explicitly, we can define the vector 
$S = {\lambda}_4 e_1 + \left(  {\lambda}_5 - \frac{1}{2} \left( 1 - m k \right) \right) e_2 - {\lambda}_4 e_3$, 
with the property that $\left( p, S \right) = 0$. 
We then define $A = \left( A_1, A_2, A_3 \right): P^{\mu} \rightarrow \mathrm{SO}(2, 1)$ by 
\begin{align*}
A_1 &= 
\frac{1}{\sqrt{2}} 
\left( \frac{p}{\Vert p \Vert} \times \frac{S}{\Vert S \Vert} + \frac{S}{\Vert S \Vert} \right),\\
A_2 &= \frac{p}{\Vert p \Vert},\\
A_3 &= 
\frac{1}{\sqrt{2}} 
\left( \frac{p}{\Vert p \Vert} \times \frac{S}{\Vert S \Vert} - \frac{S}{\Vert S \Vert} \right).
\end{align*}
Defining the map $Q: P^{\mu} \rightarrow \Mink$ by 
\[
Q = - \frac{1}{\sqrt{C_1}} \, A_2 \times v
= \frac{\langle v, A_3 \rangle}{\Vert p \Vert^2} A_1 - 
\frac{\langle v, A_1 \rangle}{\Vert p \Vert^2} A_3,
\]
we let
\[
g := (Q, A): P^{\mu} \rightarrow \widetilde{\mathcal{P}}(\mu).
\]
It then follows that $\left( p, v \right) = \mathrm{Ad}^*(g) \mu$, as required.
Therefore the map 
$\left( p, v \right) \in \widetilde{\mathcal{P}}(\mu) \rightarrow \left( g(p, v)^{-1}, (p, v) \right) \in P^{\mu}$ 
is a cross-section of the Marsden-Weinstein fibration. 

\medskip\noindent{\textbf{Negative type}}: 
For orbits of negative type we have $C_1 := \Vert p \Vert^2 < 0$. 
We treat the case where the vector $p$ is future-directed, 
although the past-directed case may be treated similarly. 
The standard form of elements in this case is  
$\mu = \left( \mathfrak{m}_1, \mathfrak{m}_2 \right)$, where 
\[
\mathfrak{m}_1 = \sqrt{\frac{|C_1|}{2}} \begin{pmatrix} 1\\ 0\\ 1 \end{pmatrix}, \quad
\mathfrak{m}_2 = - \frac{C_2}{\sqrt{2 |C_1|}} \begin{pmatrix} 1\\ 0\\ 1 \end{pmatrix}.
\]
To define a suitable basis, 
we fix a vector $S \in \Mink$ with $\Vert S \Vert^2 = 1$ and $\langle p, S \rangle = 0$. 
(For example, the vector $S =  e_2 + p^2/p^3 e_1$.) 
We then define a basis 
\[
A_1 = \frac{1}{\sqrt{2|C_1|}} \left( p - p \times S  \right), \quad
A_2 = S, \quad
A_3 = \frac{1}{\sqrt{2|C_1|}} \left( p + p \times S  \right). 
\]
Letting $A := \left( A_1, A_2, A_3 \right) \in \mathrm{SO}(2, 1)$, 
and 
\[
Q := \frac{1}{|C_1|} \, p \times v 
=  \frac{1}{\sqrt{2 |C_1|}} 
\left( \vphantom{|^|} \langle v, A_1 - A_3 \rangle A_2 - \langle v, A_2 \rangle \left( A_1 - A_3 \right) \right), 
\] 
we then define $g := (Q, A)$. 
It follows that $\eta = \mathrm{Ad}^*(g) \mu$, and hence that the map
$\left( p, v \right) \mapsto \left( g(p, v)^{-1}, (p, v) \right)$ 
is a cross-section of the Marsden-Weinstein fibration.

\medskip\noindent{\textbf{Null type}}: 
Finally, orbits of negative type have $C_1 := \Vert p \Vert^2 = 0$ with $p \neq 0$. 
Again we treat the case where the vector $p$ is future-directed, 
the past-directed case being similar. 
In this case, we may use the action of $\mathbb{E}(2, 1)$ to reduce $\mu = ( \mathfrak{m}_1, \mathfrak{m}_2)$ 
to the standard form
\[
\mathfrak{m}_1 = \begin{pmatrix} 1\\ 0\\ 0 \end{pmatrix}, \quad
\mathfrak{m}_2 = - C_2 \begin{pmatrix} 0\\ 0\\ 1 \end{pmatrix}.
\]

We now define the null basis vector $A_1 := p$, 
and extend to a basis $\left( A_1, A_2, A_3 \right)$ 
by defining, for example, 
\[
A_2 =  e_2 + \frac{p^2}{p^3} e_1, \quad 
A_3 = - \frac{1}{p^3} e_1.
\]
We let
\[
Q := A_3 \times v = \langle v, A_2 \rangle A_3 - \langle v, A_3 \rangle A_2, 
\]
and then define $g = (Q, A)$. 
Again it follows that $\eta = \mathrm{Ad}^*(g) \mu$ and therefore that the map 
$\left( p, v \right) \mapsto \left( g(p, v)^{-1}, (p, v) \right)$ 
is a cross-section of the Marsden-Weinstein fibration.

\medskip

The explicit parametrizations of the orbits given in Section~\ref{sec:parametrisation} have 
the property that $\eta^* (\sigma^{\mu}) = dt$, and hence $v^{\mu} = 1$. 
From the explicit forms of the cross-sections, $g$, it is straightforward to check for 
each type of orbit that $g(t)^{-1} g(t)^{\prime} + g(t)^{-1} \mathcal{H}[\eta(t)] g(t)$ 
lies in $\mathfrak{g}_{\mu}$ for all $t$ in the relevant range, as required. 
We may then, by direct integration, compute the gauge-transformation~\eqref{gaugetransformation}. 
The integral curves of the Euler-Lagrange system are then given by~\eqref{horizontalcurve}.

\appendix
\section{Pfaffian differential systems with one independent variable}
\begin{definition}
Let $M$ be a smooth manifold. 
A \emph{Pfaffian differential system\/} $(\mathcal{I}, \omega)$ 
with one independent variable on $M$ 
consists of a Pfaffian differential ideal ${\mathcal{I}} \subset \Omega^{*}(M)$ 
and a non-vanishing $1$-form $\omega \in \Omega^1(M)$ 
such that $\omega \nequiv 0 \pmod {\mathcal{I}}$.
\end{definition}

\begin{definition}
An \emph{integral element\/} of $(\mathcal{I}, \omega)$ is a pair $(x, E)$
consisting of a point $x \in M$ and a $1$-dimensional linear subspace
$E \subset T_x(M)$ such that $\eta \!\! \mid_{E}=0, \forall \eta \in {\mathcal{I}}$ and
$\omega \!\! \mid _{E} \neq 0$.
We denote by $V ({\mathcal{I}}, \omega)$ the set of integral elements of $({\mathcal{I}}, \omega)$.
We say that $\mathcal{I}$ has constant rank if it is generated by the cross-sections of a sub-bundle $Z$ of $T^*(M)$.
\end{definition}

\begin{definition}
A \emph{(parametrized) integral curve\/} of $(\mathcal{I}, \omega)$ is a smooth
curve $\alpha: (a, b) \subseteq \mathbb{R} \to M$ such that
\begin{gather*}
\alpha^{*}(\eta)=0, \quad \forall \eta \in \mathcal{I},\\
\gamma ^{*}(\omega) = dt.
\end{gather*}
We denote the set of integral curves of the system by $\mathcal{V}(\mathcal{I}, \omega)$.
\end{definition}

\begin{definition}
We say that the Pfaffian system in one independent variable $(\mathcal{I}, \omega)$ 
is \emph{reducible\/} if there exists a non-empty submanifold $M^* \subseteq M$
such that
\begin{itemize}
\item for each point $x \in M^*$ there exists an integral element
$(x, E) \in V(\mathcal{I}, \omega)$ tangent to $M^*$;
\item if $N \subseteq M$ is any other submanifold with the same
property then $N \subseteq M^*$.
\end{itemize}
We call $M^*$ the \emph{reduced space\/}. 
We define on $M^*$ the \emph{reduced Pfaffian system\/}, denoted by $(\mathcal{I}^*, \omega)$, 
which is obtained by restricting the original system $(\mathcal{I}, \omega)$ to $M^*$. 
\end{definition}

A basic result is the following, a proof of which may be found in \cite{Griffiths1}.
\begin{proposition}
The Pfaffian systems $(\mathcal{I}, \omega)$ and $(\mathcal{I}^*, \omega)$ have the same integral curves.
\end{proposition}

There is an algorithmic procedure for constructing the 
reduction of a Pfaffian system \cite{Griffiths1}. 
To construct the reduced space $M^*$, we consider the projection 
$M_1 \subseteq M$ of $V(\mathcal{I}, \omega)$ to $M$.
If $M_1$ is a non-empty submanifold of $M$, we then define 
$(\mathcal{I}_1, \omega_1)$ to be the restriction of $(\mathcal{I}, \omega)$ to $M_1$. 
We then construct $V(\mathcal{I}_1, \omega_1)$, the set of integral elements of $(\mathcal{I}_1, \omega_1)$. 
Repeating this construction, we inductively define
\[
M_k = (M_{k-1})_1, \quad \mathcal{I}_k = (\mathcal{I}_{k-1})_1, \quad \omega_{k} = (\omega_{k-1})_1.
\]
This process defines a sequence
$M \supseteq M_1 \supseteq \dots \supseteq M_k \supseteq \dots$ of submanifolds of $M$.
If $M^*:= \bigcap_{k \in \mathbb{N}} M_k \neq \emptyset$ then $M^*$ is the reduced space of the system.
Notice that this procedure requires that, at each stage, 
the subset $M_k \subseteq M_{k-1}$ is a non-empty submanifold.

\subsection{Cartan systems}
Let $\Psi \in \Omega^2(M)$ be an exterior differential $2$-form on $M$. 
We define the \emph{Cartan ideal\/} to be the Pfaffian differential ideal 
$\mathcal{C}(\Psi) \subseteq \Omega^*(M)$ generated by the set of 
$1$-forms $\eta_V:= i_{V} \Psi$ obtained by contracting $\Psi$ with
vector fields on $M$. If $\theta^1, \dots, \theta^n$ is a
local coframing on $M$ and if $\Psi = a_{ij} \, \theta^i \wedge \theta^j$, 
then $\mathcal{C}(\Psi)$ is locally generated by the $1$-forms $a_{ij} \, \theta^j$. 
A \emph{Cartan system\/} is a pair $(\mathcal{C}(\Psi), \omega)$ 
consisting of a Cartan ideal $\mathcal{C}(\Psi)$ and a $1$-form 
$\omega \in \Omega^1(M)$ such that $\omega|_p \notin \mathcal{C}(\Psi)|_p$, 
for all $p \in M$.

\begin{definition}
A Cartan system $(\mathcal{C}(\Psi), \omega)$ is \emph{regular\/} if
\begin{itemize}
\item it is reducible and the reduced phase space $M^*$ is odd-dimensional;
\item the $2$-form $\Psi^*:= \left. \Psi \right|_{M^*}\in \Omega^2(M^*)$ is of maximal rank on $M^*$.
\end{itemize}
\label{regular}
\end{definition}

\noindent{}An important fact is that if $(\mathcal{C}(\Psi), \omega)$ is regular 
then the Cartan ideal $\mathcal{C}(\Psi^*)$ on $M^*$ 
is the restriction to $M^*$ of the Cartan ideal $\mathcal{C}(\mathcal(\Psi))$ on $M$: 
\[
\mathcal{I}^* := \left. \mathcal{C}(\Psi)\right|_{M^*} = \mathcal{C}(\Psi^*).
\]
Again, a proof of this result may be found in \cite{Griffiths1}.

If $(\mathcal{C}(\Psi), \omega)$ is regular, 
then there exists a unique vector field $\xi$ on $M^*$ 
such that $i_{\xi} \Psi^* = 0$ and $\omega(\xi) = 1$. 
We call $\xi$ the \emph{characteristic vector field\/} of the Cartan system $(\mathcal{C}(\Psi), \omega)$.
The integral curves of the characteristic vector field coincide with the 
parametrized integral curves of $(\mathcal{C}(\Psi^*), \omega)$, 
and hence with those of $(\mathcal{C}(\Psi), \omega)$. 

\subsection{Contact systems on jet spaces}
Given a manifold $M$, we denote by $J^k(\mathbb{R}, M)$ the bundle of the $k$-order
jets of maps $\gamma: \mathbb{R} \to M$. 
The $k$-jet of $\gamma$ at $t$ will be denoted by $j^k(\gamma)|_t$. 
Local coordinates $(x^1, \dots, x^n)$ on $M$ give standard local coordinates 
$(t, x^1, \dots, x^n, x_{1}^1, \dots, x_1^n, \dots, x_k^1, \dots, x_k^n)$ 
on the jet space $J^k(\mathbb{R}, M)$. 
With respect to such a coordinate system a $k$-jet with coordinates
$(t, x^1, \dots, x^n, x_{1}^1, \dots, x_1^n, \dots, x_k^1, \dots, x_k^n)$ 
is represented by $j^k(\gamma)|_{t}$, 
where $\gamma$ is the curve defined by
\[
\gamma: s \mapsto (x^1, \dots, x^n) + (x_1^1, \dots, x_1^n)(t-s)+ \dots 
+ \frac{1}{k!}(x_k^1, \dots, x_k^n)(t-s)^k.
\]
The \emph{canonical contact system\/} $\mathcal{I}$ on $J^k(\mathbb{R}, M)$ is
defined to be the Pfaffian differential ideal generated by the forms 
$\eta^i_a = dx_a^i - x_{a+1}^i dt, i = 0, \dots, n, a = 0, \dots, k-1$
(where $x^i_0 = x^i$).
The independence condition of the system is given by the $1$-form $dt$.
The integral curves $\Gamma: (a, b) \to J^k(\mathbb{R}, M)$ of $\mathcal{I}$ such that
$\Gamma^*(\omega)=dt$ are the canonical lifts $j^k(\gamma)$ of maps $\gamma: (a, b) \to M$.

\section{Constrained variational problems in one independent variable}
\begin{definition}
Let $(\mathcal{I}, \omega)$ be a Pfaffian differential system on a smooth manifold $M$
and let $L: M \to \mathbb{R}$ be a smooth function.
The triple $(\mathcal{I}, \omega, L)$ is said to be a \emph{constrained variational problem in one independent variable\/}.
The function $L$ is referred to as the \emph{Lagrangian\/} of the variational problem.
\end{definition}

The Lagrangian $L$ gives rise to the action functional
$\mathcal{L}: \mathcal{V}(\mathcal{I}, \omega) \to \mathbb{R}$
defined (perhaps not everywhere) on the space of the integral curves of $(\mathcal{I}, \omega)$
by
\[
\mathcal{L}(\gamma) = \int_{\gamma} \gamma^* \left( L \, \omega \right).
\]

\begin{definition}
By an \emph{extremal curve\/} of $(\mathcal{I}, \omega, L)$ we mean an integral curve 
$\gamma$ that is a critical point of the functional $\mathcal{L}$ when one considers 
compactly supported variations of $\gamma$ through integral curves of the system.
\end{definition}

Let us suppose that the Pfaffian ideal $\mathcal{I}$ is generated by a sub-bundle $Z \subset T^*(M)$.
We then let ${\widetilde Z} \subset T^*(M)$ be the affine sub-bundle $L \omega + Z$.
We denote by $\psi \in \Omega^1({\widetilde Z})$ the restriction to ${\widetilde Z}$ of the tautological one-form on $T^*(M)$, 
and call $\psi$ the \emph{Liouville form\/} of the variational problem.
We let $\Psi$ be the $2$-form $d\psi$ and we consider on ${\widetilde Z}$ the Cartan system
$\mathcal{C}(\Psi)$ together with the independence condition $\omega$.

\begin{definition}
We say that $(\mathcal{I}, \omega, L)$ is a \emph{regular variational problem\/} 
if the Cartan system $(\mathcal{C}(\Psi), \omega)$ is reducible.
The reduced space $Y \subset {\widetilde Z}$ of $(\mathcal{C}(\Psi), \omega)$
is called the \emph{momentum space\/} of the variational problem.
The restriction of the Cartan system $(\mathcal{C}(\Psi), \omega)$ to $Y$ is called the
\emph{Euler-Lagrange system\/} of the variational problem, and denoted $(\mathcal{E}, \omega)$. 
\end{definition}

The importance of the Euler-Lagrange system comes from the following theorem (cf. \cite{Griffiths1,Bryant}): 

\begin{thm}
Let $\Gamma: (a, b) \to Y$ be an integral curve of the Euler-Lagrange system.
Then $\gamma = \pi_M \circ \Gamma: (a, b) \to M$ is a critical point of the action functional $\mathcal{L}$, 
where $\pi_M: Y \to M$ denotes the restriction to $Y$ of the projection $T^*(M) \to M$.
\end{thm}

This theorem allow us to find critical points of the variational problem from the integral curves the Euler-Lagrange system.
However, not all the extremals arise this way for a general variational problem.
It is known that if all the derived systems of $Z$ have constant rank,
then all the extremals are projections of the integral curves of the Euler-Lagrange system \cite{Bryant}.
Other results in this direction have been proved by L.\ Hsu \cite{Hsu}.

\begin{definition}
A variational problem $(\mathcal{I}, \omega, L)$ is said to be \emph{non-degenerate\/}
if the Cartan system $(\mathcal{C}(\Psi), \omega)$ is regular, 
in the sense of Definition~\ref{regular}.
\end{definition}

\end{document}